\documentclass{article}
\usepackage{CJK,CJKnumb,CJKulem,times,dsfont,ifthen,mathrsfs,latexsym,amsfonts,color}
\usepackage{amsmath,amsthm,makeidx,fontenc,amssymb,bm,graphicx,psfrag,listings,curves,extarrows}
\makeindex
\newtheorem{definition}{Definition}[section]
\newtheorem{theorem}{Theorem}[section]
\newtheorem{lemma}{Lemma}[section]
\newtheorem{corollary}{Corollary}[section]
\newtheorem{proposition}{Proposition}[section]
\newtheorem{remark}{Remark}[section]

\newcommand{\s}{\section}

\newcommand{\R}{\mathbb R}

\newcommand{\lab}{\label}
\newcommand{\bt}{\begin{theorem}}
\newcommand{\et}{\end{theorem}}
\newcommand{\bl}{\begin{lemma}}
\newcommand{\el}{\end{lemma}}
\newcommand{\bd}{\begin{definition}}
\newcommand{\ed}{\end{definition}}
\newcommand{\bc}{\begin{corollary}}
\newcommand{\ec}{\end{corollary}}
\newcommand{\bp}{\begin{proof}}
\newcommand{\ep}{\end{proof}}
\newcommand{\bx}{\begin{example}}
\newcommand{\ex}{\end{example}}
\newcommand{\bi}{\begin{exercise}}
\newcommand{\ei}{\end{exercise}}
\newcommand{\bo}{\begin{proposition}}
\newcommand{\eo}{\end{proposition}}
\newcommand{\br}{\begin{remark}}
\newcommand{\er}{\end{remark}}
\newcommand{\be}{\begin{equation}}
\newcommand{\ee}{\end{equation}}
\newcommand{\ba}{\begin{align}}
\newcommand{\ea}{\end{align}}
\newcommand{\bn}{\begin{enumerate}}
\newcommand{\en}{\end{enumerate}}
\newcommand{\bg}{\begin{align*}}
\newcommand{\bcs}{\begin{cases}}
\newcommand{\ecs}{\end{cases}}

\newcommand{\NN}{{\mathbb N}}

\newcommand{\bean}{\begin{eqnarray*}}
\newcommand{\eean}{\end{eqnarray*}}

\numberwithin{equation}{section}

\begin{document}
\begin{CJK*}{GBK}{song}
\title{\bf  A perturbed  nonlinear elliptic PDE  with  two Hardy-Sobolev critical exponents\thanks{Supported by NSFC (11025106, 11371212, 11271386) and the Both-Side Tsinghua Fund.
E-mails:   zhongxuexiu1989@163.com\quad \quad wzou@math.tsinghua.edu.cn}}

\date{}
\author{
{\bf X. Zhong}      \\
\footnotesize   Department of Mathematical Sciences, Tsinghua University,\\
\footnotesize  Beijing 100084, China \\\\
{\bf W. Zou}\ \\
\footnotesize Department of Mathematical Sciences, Tsinghua University,\\
\footnotesize Beijing 100084, China}

\maketitle

\vskip0.6in

\begin{center}
\begin{minipage}{120mm}
\begin{center}{\bf Abstract}\end{center}

Let  $\Omega$ be a $C^1$ open bounded domain in $\R^N$ ($N\geq 3$ )  with $0\in \partial \Omega$. Suppose that $\partial\Omega$ is $C^2$ at $0$ and the mean curvature of $\partial\Omega$ at $0$ is negative.   Consider the following perturbed PDE involving two Hardy-Sobolev critical exponents:
$$
\begin{cases}
&\Delta u+\lambda_1 \frac{u^{2^*(s_1)-1}}{|x|^{s_1}}+\lambda_2\frac{u^{2^*(s_2)-1}}{|x|^{s_2}}+\lambda_3\frac{u^p}{|x|^{s_3}}=0\;\quad \hbox{in}\;\Omega,\\
&u(x)>0\;\hbox{in}\;\Omega,\;\; u(x)=0\;\hbox{on}\;\partial\Omega,
\end{cases}
$$
where $0<s_2<s_1<2, 0\leq s_3<2, 2^*(s_i):=\frac{2(N-s_i)}{N-2}, 0\neq \lambda_i\in \R, \lambda_2>0, 1< p\leq 2^*(s_3)-1$.
The existence of ground state solution is studied under different assumptions  via  the concentration compactness principle and the  Nehari manifold method. We also apply a perturbation method to study the existence of positive solution.

 \vskip0.23in

{\it   Key  words:}  Elliptic PDE, Ground state,  Sobolev-Hardy critical exponent.\\
{\it Mathematics Subject Classification:}35J15,35J20,35J60.

\vskip0.23in

\end{minipage}
\end{center}
\vskip0.26in
\newpage
\s{Introduction}
\renewcommand{\theequation}{1.\arabic{equation}}
\renewcommand{\theremark}{1.\arabic{remark}}
\renewcommand{\thedefinition}{1.\arabic{definition}}
Consider the existence of ground state solution to the following problem
\be\lab{P1}
\begin{cases}
&\Delta u+\lambda_1 \frac{u^{2^*(s_1)-1}}{|x|^{s_1}}+\lambda_2\frac{u^{2^*(s_2)-1}}{|x|^{s_2}}+\lambda_3\frac{u^p}{|x|^{s_3}}=0\;\quad \hbox{in}\;\Omega,\\
&u(x)>0\;\hbox{in}\;\Omega,  \quad u(x)=0\;\hbox{on}\;\partial\Omega,
\end{cases}
\ee
where $ \Omega\subset \R^N$  ($N\geq 3$)   is a $C^1$ open bounded smooth domain with $0\in \partial\Omega$  and  $\partial\Omega$ is $C^2$ at $0$ and the mean curvature $H(0)<0$. The parameters satisfy  $$0<s_2<s_1<2, 0\leq s_3<2, \lambda_2>0, 1<p<2^*(s_3)-1.$$
Recall the following double critical problem
\be\lab{2014-3-25-e1}
\begin{cases}
&\Delta u+\lambda \frac{u^{2^*(s_1)-1}}{|x|^{s_1}}+\frac{u^{2^*(s_2)-1}}{|x|^{s_2}}=0\quad \hbox{in}\;\Omega,\\
&u(x)>0\;\hbox{in}\;\Omega, \quad u(x)=0\;\hbox{on}\;\partial\Omega.
\end{cases}
\ee
There has been a lot of papers concerning  (\ref{2014-3-25-e1})  under the premise of $s_2<s_1$. We note that the case of $s_1<s_2$ with $\lambda>0$ is  essentially the same.  For the case of $s_1=2 $ and $(i)N\geq 3, \lambda<(\frac{N-2}{2})^2, 0<s_2<s_1$ or $(ii) N\geq 4, 0<\lambda<(\frac{N-2}{2})^2, s_2=0$, we refer to \cite{ChernLin.2010,GhoussoubKang.2004,GhoussoubRobert.2006}.
When $s_2=0$, equation (\ref{2014-3-25-e1}) becomes
 \be\lab{2014-3-10-e1}
 \Delta u+\lambda \frac{u^{2^*(s_1)-1}}{|x|^{s_1}}+u^{\frac{N+2}{N-2}}=0\;\hbox{in}\;\Omega.
 \ee
It is well known that (\ref{2014-3-10-e1}) has no least-energy solution  if $0\leq s_1<2$ with $\lambda<0$. However, for the case of
$\lambda>0, 0<s_1<2$ and $s_2=0$, the existence of positive solution is proved by Hsia, Lin  and Wadade  \cite{HsiaLinWadade.2010}.
In the very  recent paper \cite{LiLin.2012}, the existence of positive solution for $N\geq 3,\lambda\in \R, 0<s_2<s_1<2$ is proved by Li  and Lin.
Basically, (\ref{2014-3-25-e1}) has been studied for all the choices of the parameters $s_1,s_2$ under the premise that the coefficient of the highest  power  term is positive.  However,  an
open problem is proposed   by  Li and Lin  in  \cite[Remark 1.2]{LiLin.2012}  which says:   {\it  For the situation
$s_1<s_2$ and $\lambda<0$, the existence of positive solutions to (\ref{2014-3-25-e1}) is completely open. Even for the equation
\be\lab{zou=1} \Delta u-u^p+\frac{u^{2^*(s)-1}}{|x|^s}=0 \;\hbox{in}\;    \Omega, \ee
where $0<s<2$ and $2^*(s)-1<p<\frac{N+2}{N-2}$, the existence problem still remains an interesting open question. }  It seems that the first partial answer to this open problem
is obtained in   \cite{CeramiZhongZou.2014}.

\vskip0.12in

Further, although (\ref{2014-3-10-e1}) has no least-energy solutions for $\lambda<0, 0<s_1<2$,  the following   perturbed equation
\be\lab{2014-3-25-e2}
\begin{cases}
&\Delta u-\frac{u^{2^*(s)-1}}{|x|^s}+u^{\frac{N+2}{N-2}}+u^p=0\;\hbox{in}\;\Omega,\\
&u(x)>0\;\hbox{in}\;\Omega\;\hbox{and}\;u(x)=0\;\hbox{on}\;\partial\Omega
\end{cases}
\ee
has a positive solution if $N\geq 4, 2^*(s)-1<p<\frac{N+2}{N-2}$, see   Li and Lin  \cite[Theorem 5.1]{LiLin.2012}.



\vskip0.11in

In the current paper, we are interested in the more general perturbation problem than   (\ref{2014-3-25-e2}), that is, the equation   (\ref{P1}).   We obtain the following main theorems:

\bt\lab{main-th}
Suppose  that     $\Omega$    is an open bounded smooth domain in $\R^N$  $(N\geq 3)$, $0\in \partial\Omega$ and the mean curvature  of  $\partial\Omega$  at $0$ is negative,  i.e.,    $H(0)<0$. Assume that $  0<s_2<s_1<2, 0\leq s_3<2, \lambda_2>0, 1<p<2^*(s_3)-1$, and  that one of the following conditions is satisfied:
\begin{itemize}
\item[(1)]$\lambda_1>0, \lambda_3>0$.
\item[(2)] $\lambda_1>0, \lambda_3<0, p\leq 2^*(s_1)-1$.
\item[(3)] $\lambda_1<0, \lambda_3>0, p\geq 2^*(s_1)-1$.
\item[(4)] $\lambda_1<0, \lambda_3<0, p<2^*(s_2)-1$.
\end{itemize}
Furthermore, if $\lambda_3<0$, we require  either $p<\frac{N-s_3}{N-2}$ or $p\geq\frac{N-s_3}{N-2}$ with $|\lambda_3|$ small enough. Then (\ref{P1}) possesses a ground state solution.
\et

\br
 We remark  that Theorem \ref{main-th} does not cover the  following two cases:
 \begin{itemize}
 \item   $\lambda_3<0,$   $2^*(s_1)-1<p\leq 2^*(s_3)-1$  and  $\lambda_1>0$;
  \item   $\lambda_3<0,$  $2^*(s_2)-1<p\leq 2^*(s_3)-1$   and $\lambda_1<0$.
  \end{itemize}
 Since for these cases,    we do not know whether   the  $(PS)$ sequence is  bounded or not.  In particular, the Nehari manifold method fails.   The existence of  the ground state solution  for this two  cases  remains open. \er

However,  when  $|\lambda_3|$ small enough  we may  obtain the existence of positive solution. Precisely, we have the following result:
\bt\lab{main-th2}
Suppose that  $\Omega$ is an open  bounded smooth domain in $\R^N$  ($N\geq 3)$,  $0\in \partial\Omega$ and the mean curvature $H(0)<0$. Assume that $$ 0<s_2<s_1<2, 0\leq s_3<2, \lambda_2>0, \lambda_3<0$$
and that either  $2^*(s_1)-1<p\leq 2^*(s_3)-1$ if $\lambda_1>0$ or $2^*(s_2)-1<p\leq 2^*(s_3)-1$ when  $\lambda_1<0$.  Then there exists $\lambda_0<0$  such that (\ref{P1}) has a positive solution for $\lambda_0<\lambda_3<0$.
\et

\br   In the above Theorem   \ref{main-th2}, we allow     $p= 2^*(s_3)-1$, it means that the equation  (\ref{P1}) has three Hardy-Sobolev critical terms.\er

This paper is organized as follows. In section 2, we give some properties of the Nehari manifold. Since the problem involves critical terms, it is well known that the lack of the compactness will bring much troubles.
In section 3, we will  determine the threshold of
the functional for which the Palais-Smale condition holds and check that the ground state value lies in the safe region.
Based  on these  preparations, we prove Theorem \ref{main-th}. In section 5, we will prove Theorem \ref{main-th2} by a perturbation method.



\s{Nehari manifold}
\renewcommand{\theequation}{2.\arabic{equation}}
\renewcommand{\theremark}{2.\arabic{remark}}
\renewcommand{\thedefinition}{2.\arabic{definition}}
Let $L^{p}(\Omega,\frac{dx}{|x|^s})$ denote the space of $L^{p}$-integrable functions with respect to the measure $\frac{dx}{|x|^s}$.   Let $|u|_{s, p}:=\big(\int_\Omega \frac{|u|^p}{|x|^s}dx\big)^{\frac{1}{p}}$ and $|u|_{p}:=|u|_{0,p}$. The Hardy-Sobolev inequality  (see  \cite{CaffarelliKohnNirenberg.1984,CatrinaWang.2001,GhoussoubYuan.2000})
asserts that $D_{0}^{1,2}(\R^N)\hookrightarrow L^{2^*(s)}(\R^N,\frac{dx}{|x|^s})$ is a  continuous embedding for $s\in [0,2]$. That is, there exists $C_s>0$ such that
\be\lab{2014-2-27-e1}
\Big(\int_{\R^N}\frac{|u|^{2^*(s)}}{|x|^s}dx\Big)^{\frac{2}{2^*(s)}}\leq C_s \int_{\R^N}|\nabla u|^2dx\;\hbox{for all}\;u\in D_{0}^{1,2}(\R^N).
\ee
If $|\Omega|<\infty$    and $p<2^*(s_3)-1$, we can obtain that
$\displaystyle \int_\Omega \frac{|u|^{p+1}}{|x|^{s_3}}dx<\infty$ for all $u\in H_0^1(\Omega)$. In particular, the embedding $H_0^1(\Omega)\hookrightarrow L^{p+1}(\Omega,\frac{dx}{|x|^{s_3}})$ is compact which was  established in \cite[Theorem 1.9]{Willem.1996}  for $s_3=0$ and \cite[Lemma 2.1]{CeramiZhongZou.2014} for $0<s_3<2$.  A function $u\in H_0^1(\Omega)$ is said to be a weak solution  to the problem (\ref{P1}) iff
$$\int_\Omega \nabla u\cdot \nabla v dx-\lambda_1\int_\Omega \frac{|u|^{2^*(s_1)-2}uv}{|x|^{s_1}}dx-$$
\be\lab{2014-3-18-e1}
 \lambda_2\int_\Omega \frac{|u|^{2^*(s_2)-2}uv}{|x|^{s_2}}dx-\lambda_3\int_\Omega \frac{|u|^{p-1}uv}{|x|^{s_3}}dx=0
\ee
for all $v\in H_0^1(\Omega)$.
Thus, the corresponding energy functional of (\ref{P1}) is
\be\lab{2014-3-18-e2}
\Phi(u)=\frac{1}{2}\|u\|^2-\frac{\lambda_1}{2^*(s_1)}|u|_{{s_1},2^*(s_1)}^{2^*(s_1)}
-\frac{\lambda_2}{2^*(s_2)}|u|_{{s_2},2^*(s_2)}^{2^*(s_2)}
-\frac{\lambda_3}{p+1}|u|_{{s_3},p+1}^{p+1}.
\ee
The associated Nehari manifold is defined as
$$\mathcal{N}:=\Big\{u\in H_0^1(\Omega)\backslash \{0\}:J(u)=0\Big\},$$
where \begin{align}\lab{2014-3-18-e3}
J(u):=\langle \Phi'(u), u\rangle
=\|u\|^2-\lambda_1|u|_{{s_1},2^*(s_1)}^{2^*(s_1)}-\lambda_2|u|_{{s_2},2^*(s_2)}^{2^*(s_2)}
-\lambda_3|u|_{{s_3},p+1}^{p+1}
\end{align}
and $\Phi'(u)$ denotes the Fr\'{e}chet derivative of $\Phi$ at $u$;  $\langle \cdot, \cdot\rangle$ is the dual  product between $H_0^1(\Omega)$ and its dual space $H^{-1}(\Omega)$. We have the following properties on the Nehari manifold.

\bl\lab{Nehari-l1}
Assume that $0<s_2<s_1<2, 0\leq s_3<2, \lambda_2>0, 1<p< 2^*(s_3)-1$. Then $\forall\;u\in H_0^1(\Omega)\backslash\{0\}$, there exists a unique $t=t_{u}>0$ such that $tu\in \mathcal{N}$ if one of the following assumptions is satisfied:
\begin{itemize}
\item[(1)]$\lambda_1>0, \lambda_3>0$.
\item[(2)] $\lambda_1>0, \lambda_3<0, p\leq 2^*(s_1)-1$.
\item[(3)] $\lambda_1<0, \lambda_3>0, p\geq 2^*(s_1)-1$.
\item[(4)] $\lambda_1<0, \lambda_3<0, p<2^*(s_2)-1$.
\end{itemize}
Moreover, $\mathcal{N}$ is closed and bounded away from $0$.
\el
\bp
For any $u\in H_0^1(\Omega)$, we denote
\be\lab{2014-3-19-e1}
a(u):=\|u\|^2,\;
b(u):=|u|_{{s_1},2^*(s_1)}^{2^*(s_1)},\;
c(u):=|u|_{{s_2},2^*(s_2)}^{2^*(s_2)},\;
d(u):=|u|_{{s_3},p+1}^{p+1}.
\ee
We will write them as $a, b, c, d$ for simplicity if there is no ambiguity.
Then, $\displaystyle\frac{d}{dt}\Phi(tu)=tg(t)$, where
$$g(t):=a-\lambda_1 b t^{2^*(s_1)-2}-\lambda_2 ct^{2^*(s_2)-2}-\lambda_3dt^{p-1}.$$    We also see that for $t>0, \frac{d}{dt}\Phi(tu)=0$ if and only if $g(t)=0$.
Recalling that $\lambda_2>0, s_2<s_1$ and $ p<2^*(s_2)-1$ if $\lambda_3<0$, we obtain that $g(t)\rightarrow -\infty$ as $t\rightarrow +\infty$. Combine  with $g(0)=a>0$, we have that there exists some $t>0$ such that $g(t)=0$ due to the continuity of $g(t)$.   It follows that $tu\in \mathcal{N}$.  Let $u\in \mathcal{N}$, since $p>1, 2^*(s_i)>2$, by the embedding theorem we obtain that
\begin{align*}
a=\lambda_1 b+\lambda_2c+\lambda_3d
\leq C \Big(a^{\frac{2^*(s_1)}{2}}+a^{\frac{2^*(s_2)}{2}}+a^{\frac{p+1}{2}}\Big),
\end{align*}
which implies that there exists some $\delta_0>0$ such that
\be\lab{2014-3-18-e4}
\|u\|=a^{\frac{1}{2}}\geq \delta_0\;\hbox{for all}\;u\in \mathcal{N}.
\ee
Then for any $u\neq 0$, $t_0:=\inf\{t|g(t)=0\}>0$ and by the continuity of $g(t)$, we obtain that $g(t_0)=0$. Without loss of generality, we may assume that $t_0=1$, that is, $g(1)=0$ and $g(t)>0$ for all $t\in (0,1)$.

\vskip0.1in
  If $\lambda_1>0, \lambda_3>0$, it is easy to see that
$g'(t)<0.$
Hence, $g(t)<g(1)=0$ for all $t>1$.

    If $\lambda_1>0, \lambda_3<0, p\leq 2^*(s_1)-1$, we consider $t>1$ first.  We have
$g'(t)=-t^{p-2}h(t),$
where
$$h(t):=\lambda_1 b \big(2^*(s_1)-2\big)t^{2^*(s_1)-p-1}+\lambda_2c\big(2^*(s_2)-2\big)t^{2^*(s_2)-p-1}+\lambda_3d(p-1).$$
Recall  that
$a-\lambda_1b-\lambda_2c-\lambda_3d=0,$
we obtain that
\begin{align}\lab{2014-3-18-e6}
h(t)>&\lambda_1 b \big(2^*(s_1)-2\big)+\lambda_2c\big(2^*(s_2)-2\big)+\lambda_3d(p-1)\nonumber\\
=&\lambda_1b\big(2^*(s_1)-p-1\big)+\lambda_2c\big(2^*(s_2)-p-1\big)+a(p-1)>0.
\end{align}
It follows that $g'(t)<0$ for all $t>1$ and then   $g(t)<g(1)=0\;\hbox{for all}\; t>1.$

\vskip0.1in

  If $\lambda_1<0, \lambda_3>0, 2^*(s_1)-1\leq p$, then we have
$
g'(t)=-t^{2^*(s_1)-3}q(t),
$
where
$$q(t):=\lambda_1\big(2^*(s_1)-2\big)b+\lambda_2\big(2^*(s_2)-2\big)ct^{2^*(s_2)-2^*(s_1)}+
\lambda_3(p-1)dt^{p+1-2^*(s_1)}.$$
Assume  $t>1$, we obtain that
\begin{align}\lab{2014-3-18-e8}
q(t)>&\lambda_1\big(2^*(s_1)-2\big)b+\lambda_2\big(2^*(s_2)-2\big)c+\lambda_3(p-1)d\nonumber\\
=&a\big(2^*(s_1)-2\big)+\lambda_2c\big(2^*(s_2)-2^*(s_1)\big)+\lambda_3d\big(p+1-2^*(s_1)\big)\nonumber\\
>&0.
\end{align}
Hence, we also obtain that $g'(t)<0$ for $t>1$.

\vskip0.1in

   If $\lambda_1<0,\lambda_3<0, p<2^*(s_2)-1$, we have
$$tg'(t)=-\lambda_1b\big(2^*(s_1)-2\big)t^{2^*(s_1)-2}-\lambda_3d(p-1)t^{p-1}-\lambda_2c\big(2^*(s_2)-2\big)t^{2^*(s_2)-2}.$$
Recalling that  $\lambda_2c=a-\lambda_1b-\lambda_3d,$
we obtain that
\begin{align*}
tg'(t)=&\big[\big(2^*(s_2)-2\big)t^{2^*(s_2)-2^*(s_1)}-\big(2^*(s_1)-2\big)\big]\lambda_1 b t^{2^*(s_1)-2}\\
&+\big[\big(2^*(s_2)-2\big)t^{2^*(s_2)-p-1}-(p-1)\big]\lambda_3dt^{p-1}\\
&-a\big(2^*(s_2)-2\big)t^{2^*(s_2)-2}\\
=:&I+II+III.
\end{align*}
Since $s_2<s_1$, we have $2^*(s_2)>2^*(s_1)$. Thus, for $t>1$, we obtain $$\big(2^*(s_2)-2\big)t^{2^*(s_2)-2^*(s_1)}-\big(2^*(s_1)-2\big)>2^*(s_2)-2^*(s_1)>0.$$
It follows that $I<0$ due to the fact of $\lambda_1<0$.
Similarly, since $p<2^*(s_2)-1$, we can prove that $II<0$ for $t>1$. Obviously,  $III<0$. We deduce that
$g'(t)<0\;\hbox{for all}\;t>1.$
Based  on the  above arguments, we obtain that $g(t)<g(1)=0$ for all $t>1$. Hence, for any $0\neq u\in H_0^1(\Omega)$, there exists a  unique $t>0$ denoted by $t_{u}$ such that $t_{u}u\in \mathcal{N}$. By (\ref{2014-3-18-e4}), we have that $\mathcal{N}$ is bounded away from $0$ and that $\mathcal{N}$ is closed.
\ep

\bl\lab{2014-3-25-l1}
Under the assumptions of Lemma \ref{Nehari-l1},   any $(PS)_c$ sequence  $\{u_n\}$  of $\Phi(u)$,  i.e.,
$\Phi(u_n)\rightarrow c,\;\;  \Phi'(u_n)\rightarrow 0\;\hbox{in}\;H^{-1}(\Omega) ,$
is bounded in $H_0^1(\Omega)$.
\el
\bp
Let $\{u_n\}\subset H_0^1(\Omega)$ be a $(PS)_c$ sequence of $\Phi(u)$, then we have
\be\lab{2014-3-25-x-e1}
\Phi(u_n)=\frac{1}{2}a(u_n)-\frac{\lambda_1}{2^*(s_1)}b(u_n)
-\frac{\lambda_2}{2^*(s_2)}c(u_n)-\frac{\lambda_3}{p+1}d(u_n)=c+o(1)
\ee
and
$
\langle\Phi'(u_n),u_n\rangle =a(u_n)-\lambda_1 b(u_n)-\lambda_2c(u_n)-\lambda_3 d(u_n)=o(1)\|u_n\|,$
where $a(u), b(u), c(u), d(u)$ are defined by (\ref{2014-3-19-e1}).

\vskip0.2in

\noindent {\it  Case 1.}    Assume $\lambda_1>0, \lambda_3>0$.    If $p+1\geq 2^*(s_1)$, we have that
\begin{align*}
c+o(1)(1+\|u_n\|)=\Phi(u_n)-\frac{1}{2^*(s_1)}\langle \Phi'(u_n), u_n\rangle
\geq(\frac{1}{2}-\frac{1}{2^*(s_1)})\|u_n\|^2.
\end{align*}
If $p+1<2^*(s_1)$,  note that  $s_2<s_1, p+1<2^*(s_1)$,
it follows   that
\begin{align*}
c+o(1)(1+\|u_n\|)=\Phi(u_n)-\frac{1}{p+1}\langle \Phi'(u_n), u_n\rangle
\geq(\frac{1}{2}-\frac{1}{p+1})\|u_n\|^2,
\end{align*}

\noindent {\it  Case 2.} If $\lambda_1>0, \lambda_3<0, p\leq 2^*(s_1)-1$, we have
\begin{align*}
c+o(1)=\Phi(u_n)
 =(\frac{1}{2}-\frac{1}{p+1})a(u_n)+(\frac{1}{p+1}-\frac{1}{2^*(s_1)})\lambda_1b(u_n)\\
+(\frac{1}{p+1}-\frac{1}{2^*(s_2)})\lambda_2c(u_n)+o(1)\|u_n\|.
\end{align*}
Hence,
$\displaystyle c+o(1)(1+\|u_n\|)\geq (\frac{1}{2}-\frac{1}{p+1})\|u_n\|^2.$

\noindent {\it  Case 3.}  If $\lambda_1<0, \lambda_3>0, 2^*(s_1)-1\leq p$, then
\begin{align*}
c+o(1)=&\Phi(u_n)
=\big(\frac{1}{2}-\frac{1}{2^*(s_1)}\big)a(u_n)+\big(\frac{1}{2^*(s_1)}-\frac{1}{2^*(s_2)}\big)\lambda_2c(u_n)\\
&+\big(\frac{1}{2^*(s_1)}-\frac{1}{p+1}\big)\lambda_3d(u_n)+o(1)\|u_n\|.
\end{align*}
Hence,
$\displaystyle c+o(1)(1+\|u_n\|)\geq \big(\frac{1}{2}-\frac{1}{2^*(s_1)}\big)\|u_n\|^2.$

\noindent {\it Case 4.}  If $\lambda_1<0,\lambda_3<0, p<2^*(s_2)-1$,  similarly  we have
$$c+o(1)(1+\|u_n\|)\geq \big(\frac{1}{2}-\frac{1}{2^*(s_2)}\big)\|u_n\|^2.$$
Based  on the  above  arguments, we can see that $\{u_n\}$ is bounded in $H_0^1(\Omega)$.
\ep

\br\lab{2014-3-18-r1}

Under the assumptions of Lemma \ref{Nehari-l1},  we define
\be\lab{2014-3-18-e9}
c_0:=\inf_{u\in \mathcal{N}}\Phi(u)
\ee
and
$$\eta:=\min\{\frac{1}{2}-\frac{1}{p+1}, \frac{1}{2}-\frac{1}{2^*(s_1)}\}>0.$$
Similar to the prove of Lemma \ref{2014-3-25-l1}, we see that
$c_0\geq \eta \delta_0^2>0,$
where $\delta_0$ is given by (\ref{2014-3-18-e4}).
If $c_0$ is achieved by some $u\in \mathcal{N}$, then $u$ is a ground state solution of (\ref{P1}).
\er
\bl\lab{l3}
Under the assumptions of Lemma \ref{Nehari-l1}, let $\{u_n\}\subset \mathcal{N}$ be a $(PS)_c$ sequence for $\Phi\big|_{\mathcal{N}}$, that is, $\Phi(u_n)\rightarrow c$ and $\Phi'\big|_{\mathcal{N}}(u_n)\rightarrow 0$ in $H^{-1}(\Omega)$. Then $\{u_n\}$ is also a $(PS)_c$  sequence for $\Phi$.
\el
\bp
Firstly,  by the similar arguments  as  that in Lemma \ref{2014-3-25-l1}, we may show  that $\{(u_n,v_n)\}$ is bounded in $\mathscr{D}$.
Let $\{t_n\}\subset \R$ be a sequence of multipliers satisfying
$$\Phi'(u_n)=\Phi'\big|_{\mathcal{N}}(u_n)+t_nJ'(u_n).$$
Testing by $u_n$, we obtain that $t_n \langle J'(u_n), u_n\rangle \rightarrow 0$.
Recalling that for any $u\in \mathcal{N}$, we have
$$\langle J'(u), u\rangle=2a-2^*(s_1)\lambda_1b-2^*(s_2)\lambda_2c-(p+1)\lambda_3d$$
and
$a-\lambda_1b-\lambda_2c-\lambda_3d=0,$
where $a,b,c,d$ is defined by (\ref{2014-3-19-e1}).

\noindent(1)  If $\lambda_1>0, \lambda_3>0$,  then
\begin{align*}
\langle J'(u), u\rangle=&2a-2^*(s_1)\lambda_1b-2^*(s_2)\lambda_2c-(p+1)\lambda_3d\\
<&-\min\{2^*(s_1)-2, 2^*(s_2)-2, p-1\}(\lambda_1b+\lambda_2c+\lambda_3d)\\
=&-\min\{2^*(s_1)-2, 2^*(s_2)-2, p-1\}a.
\end{align*}

\noindent{(2)} If $\lambda_1>0, \lambda_3<0, p\leq 2^*(s_1)-1$, we have $p<2^*(s_2)-1$ and then
\begin{align*}
\langle J'(u), u\rangle=2a-2^*(s_1)\lambda_1b-2^*(s_2)\lambda_2c-(p+1)(a-\lambda_1b-\lambda_2c)
<-(p-1)a.
\end{align*}

\noindent{(3)}  If $\lambda_1<0, \lambda_3>0, 2^*(s_1)-1\leq p$,
we also have
\begin{align*}
\langle J'(u), u\rangle
=&\big(2-2^*(s_1)\big)a+\big(2^*(s_1)-2^*(s_2)\big)\lambda_2c+\big(2^*(s_1)-p-1\big)\lambda_3d\\
<&-\big(2^*(s_1)-2\big)a.
\end{align*}

\noindent{(4)} If $\lambda_1<0,\lambda_3<0, p<2^*(s_2)-1$, we have
\begin{align*}
\langle J'(u), u\rangle
=&\big(2-2^*(s_2)\big)a+\big(2^*(s_2)-2^*(s_1)\big)\lambda_1b+\big(2^*(s_2)-p-1\big)\lambda_3d\\
<&-\big(2^*(s_2)-2\big)a.
\end{align*}
Thus, under the assumptions of Lemma \ref{Nehari-l1},
 $\langle J'(u), u\rangle <-\varrho \|u\|^2\;\hbox{for all}\;u\in \mathcal{N},$
 where $\varrho:=\min\{2^*(s_1)-2, 2^*(s_2)-2, p-1\}>0$.  Invoke  (\ref{2014-3-18-e4}), we have
 $$\langle J'(u_n), u_n\rangle <-\varrho \delta_0^2\;\hbox{for all}\;n.$$
Hence, we obtain that $t_n$ is bounded.
On the other hand, it is easy to see that $\langle J'(u_n), u_n\rangle$ is bounded due to the boundedness of $\{u_n\}$.
We claim that $t_n\rightarrow 0$. If not,
up to a subsequence, we may assume that $t_n\rightarrow t_0\neq 0$ and $\langle J'(u_n), u_n\rangle\rightarrow d_0<-\varrho \delta_0^2$. Then $$\big|t_n \langle J'(u_n), u_n\rangle\big|\rightarrow |t_0 d_0|>|t_0|\varrho \delta_0^2\neq 0,$$
a contradiction.
Thus, we see  that $t_n\rightarrow 0$ and it follows that $\Phi'(u_n)\rightarrow 0$ in $H^{-1}(\Omega)$.
\ep


\s{Analysis of the  Palais-Smale sequences}
\renewcommand{\theequation}{3.\arabic{equation}}
\renewcommand{\theremark}{3.\arabic{remark}}
\renewcommand{\thedefinition}{3.\arabic{definition}}
Understanding asymptotic behavior is usually  fundamental in the resolution of mathematical
problems, particularly the problem possesses critical terms.
The following result is due to \cite{LiLin.2012}:\\

\vskip 0.06in
\noindent{\bf Theorem A. }(\cite[Theorem 1.2]{LiLin.2012}) {\it  Let $N\geq 3, 0<s_2<s_1<2, \lambda\in \R$, then the following problem
\be\lab{2014-3-19-e2}
\begin{cases}
\Delta u+\lambda \frac{u^{2^*(s_1)-1}}{|x|^{s_1}}+\frac{u^{2^*(s_2)-1}}{|x|^{s_2}}=0\quad &\hbox{in}\;\R_+^N,\\
u(x)>0\quad \hbox{in}\;\Omega, \quad u(x)=0\quad &\hbox{on}\;\partial\R_+^N,
\end{cases}
\ee
has a least-energy solution $u\in H_0^1(\R_+^N)$.}\hfill$\Box$

\vskip0.1in

Let $u>0$ be  the least energy solution of (\ref{2014-3-19-e2}), then
\be\lab{2014-3-19-e3}
|u(y)|\leq C |y|^{1-N}\quad \hbox{for}\;|y|\geq 1
\hbox{ and  }
|\nabla u(y)|\leq |y|^{-N}\quad \hbox{for}\;|y|\geq 1.
\ee
See  \cite[Page 16-17]{LiLin.2012}.
We also note that,  by the well-known   moving plane method, one can prove that $u(x', x_N)$ is axially symmetric with respect to the $x_N$-axis, i.e., $u(x', x_N)=u(|x'|, x_N)$, where $x'=(x_1,\cdots,x_{N-1})$. Since the argument is standard, we omit the proof here (see  \cite[Lemma 2.6]{LinWadadeothers.2012}).
When $\lambda_2>0$, define $v=\lambda_{2}^{\frac{1}{2^*(s_2)-2}}u$, a direct calculation shows that $u$ is a solution of
\be\lab{2014-3-19-e5}
\begin{cases}
\Delta u+\lambda_1 \frac{u^{2^*(s_1)-1}}{|x|^{s_1}}+\lambda_2\frac{u^{2^*(s_2)-1}}{|x|^{s_2}}=0\quad &\hbox{in}\;\R_+^N,\\
u(x)>0\quad \hbox{in}\;\R_+^N,\quad u(x)=0\quad &\hbox{on}\;\partial\R_+^N.
\end{cases}
\ee
if and only if $v$ is a solution to (\ref{2014-3-19-e2}) with $\lambda=\lambda_1 \lambda_{2}^{\frac{2-2^*(s_1)}{2^*(s_2)-2}}$.
We denote the least energy corresponding to (\ref{2014-3-19-e5}) by $c_{\lambda_1,\lambda_2}$, that is,
$$c_{\lambda_1,\lambda_2}=\inf\{A_{\lambda_1,\lambda_2}(u)|u\;\hbox{is a solution to}\;(\ref{2014-3-19-e5})\},$$
where
$$A_{\lambda_1,\lambda_2}(u):=\frac{1}{2}\int_{\R_+^N}|\nabla u|^2dx-\frac{\lambda_1}{2^*(s_1)}\int_{\R_+^N}\frac{|u|^{2^*(s_1)}}{|x|^{s_1}}dx
-\frac{\lambda_2}{2^*(s_2)}\int_{\R_+^N}\frac{|u|^{2^*(s_2)}}{|x|^{s_2}}dx.$$
It is easy to see that
\be\lab{2014-3-19-e6}
A_{\lambda,1}(v)=\lambda_{2}^{\frac{2}{2^*(s_2)-2}} A_{\lambda_1,\lambda_2}(u),
\ee
where
$$v=\lambda_{2}^{\frac{1}{2^*(s_2)-2}}u, \lambda=\lambda_1 \lambda_{2}^{\frac{2-2^*(s_1)}{2^*(s_2)-2}}.$$
It follows that
\be\lab{2014-3-19-e7}
c_{\lambda_1,\lambda_2}=\lambda_{2}^{\frac{-2}{2^*(s_2)-2}}c_{\lambda,1}, \quad \lambda=\lambda_1 \lambda_{2}^{\frac{2-2^*(s_1)}{2^*(s_2)-2}}.
\ee
Let $w>0$   be   a ground state solution to (\ref{2014-3-19-e5}), then
\be\lab{2014-3-19-e8}
|w(y)|\leq C |y|^{1-N}\quad \hbox{for}\;|y|\geq 1
\ee
and
\be\lab{2014-3-19-e9}
|\nabla w(y)|\leq \lambda_{2}^{-\frac{1}{2^*(s_2)-2}}|y|^{-N}\quad \hbox{for}\;|y|\geq 1.
\ee
Similar to \cite[Theorem 3.1]{CeramiZhongZou.2014}, we can establish the following splitting result which provides a precise description of a behavior of $(PS)_c$ sequence for $\Phi(u)$.

\bt\lab{2014-3-19-compactness-th1}(Splitting  Theorem)
Suppose   that $\{u_n\}\subset H_0^1(\Omega)$ is a bounded $(PS)_c$ sequence of the functional $\Phi(u)$. That is, $\Phi(u_n)\rightarrow c$ and $\Phi'(u_n)\rightarrow 0$ strongly in $H^{-1}(\Omega)$ as $n\rightarrow \infty$. Then there exists a solution $U^0$ to the equation in (\ref{P1})  ($U^0\equiv 0$ is allowed), number $k\in \NN\cup \{0\},$  $k$ functions $U^1, \cdots, U^k$ and $k$ sequences of radius $r_{n}^{j}>0, 1\leq j\leq k$ such that the following properties are satisfied up  to a subsequence if necessary:  Either
\begin{itemize}
\item[(a)] $u_n\rightarrow U^0$ in $H_0^1(\Omega)$ or
\item[(b)] the following items all are true:
\begin{itemize}
\item[(b1)]$U^j\in D^{1,2}(\R_+^N)\subset D^{1,2}(\R^N)$ are nontrivial solutions of (\ref{2014-3-19-e5});
\item[(b2)]$r_{n}^{j}\rightarrow 0$ as $n\rightarrow \infty$;
\item[(b3)] $\|u_n-U^0-\sum_{j=1}^{k} (r_{n}^{j})^{\frac{2-N}{2}}U^{j}(\frac{\cdot}{r_{n}^{j}})\|\rightarrow 0$, where $\|\cdot\|$ is the norm in $D^{1,2}(\R^N)$;
\item[(b4)]$\|u_n\|^2\rightarrow \|U^0\|^2+ \sum_{j=1}^{k}\|U^{j}\|^2$;
\item[(b5)] $\Phi(u_n)\rightarrow \Phi(U^0)+\sum_{j=1}^{k}A_{\lambda_1,\lambda_2}(U^j)$,
    where
    $$A_{\lambda_1,\lambda_2}(u):=\frac{1}{2}\int_{\R_+^N}|\nabla u|^2dx-\frac{\lambda_1}{2^*(s_1)}\int_{\R_+^N}\frac{|u|^{2^*(s_1)}}{|x|^{s_1}}dx$$
    $$
-\frac{\lambda_2}{2^*(s_2)}\int_{\R_+^N}\frac{|u|^{2^*(s_2)}}{|x|^{s_2}}dx.$$
\end{itemize}
\end{itemize}
\et\hfill$\Box$

  The following corollary is a straightforward consequence  of  the above theorem.
\bc\lab{2014-3-19-cro1}
Under the assumptions of Lemma \ref{Nehari-l1},   the functional
$\Phi(u)$ satisfies $(PS)_c$ condition for $c<c_{\lambda_1,\lambda_2}$.
\ec
\bp
Let $\{u_n\}\subset H_0^1(\Omega)$ be such that $\Phi(u_n)\rightarrow c<c_{\lambda_1,\lambda_2}$, $\Phi'(u_n)\rightarrow 0$ in $H^{-1}(\Omega)$. By Lemma  \ref{2014-3-25-l1}, we obtain that $\{u_n\}$ is bounded in $H_0^1(\Omega)$.  By Theorem \ref{2014-3-19-compactness-th1}, we obtain that $u_n\rightarrow U^0$ in $H_0^1(\Omega)$ up to a subsequence. If not, $k\neq 0$  and
$$\sum_{j=1}^{k}A_{\lambda_1,\lambda_2}(U^j)
    \geq c_{\lambda_1,\lambda_2}.$$ Recalling that $c_0>0$ (see Remark \ref{2014-3-18-r1}), we have $\Phi(U^0)\geq 0$, then
$$
c=\Phi(u_n)+o(1)=\Phi(U^0)+\sum_{j=1}^{k}A_{\lambda_1,\lambda_2}(U^j)\geq
 c_{\lambda_1,\lambda_2}, $$
a contradiction.
\ep

 Before giving the proof of   Theorem  \ref{2014-3-19-compactness-th1}, we need to prepare the following two auxiliary results:

\bl\lab{2014-3-19-l2}
Let  $\{u_n\}\subset H_0^1(\Omega)$ be such that $A_{\lambda_1,\lambda_2}(u_n)\rightarrow c$ and $A'_{\lambda_1,\lambda_2}(u_n)\rightarrow 0$ in $H^{-1}(\Omega)$.
For $\{r_n\}\subset (0,\infty)$ with $r_n\rightarrow 0$,  let
$v_n(x):=r_{n}^{\frac{N-2}{2}}u_n(r_nx)$ be such that $v_n\rightharpoonup v$ in $D^{1,2}(\R^N)$ and $v_n\rightarrow v$ a.e. on $\R^N$. Then, $A'_{\lambda_1,\lambda_2}(v)=0$ and the sequence
$$w_n(x):=u_n(x)-r_{n}^{\frac{2-N}{2}}v(\frac{x}{r_n})$$
satisfies $A_{\lambda_1,\lambda_2}(w_n)\rightarrow c-A_{\lambda_1,\lambda_2}(v), A'_{\lambda_1,\lambda_2}(w_n)\rightarrow 0$ in $H^{-1}(\Omega)$ and $\|w_n\|^2=\|u_n\|^2-\|v\|^2+o(1)$.
\el
\bp
Without loss of generality, we assume  that $\partial R_+^N:=\{x_N=0\}$ is tangent to $\partial\Omega$ at $0$, and that $-e_N=(0,\cdots, -1)$ is the outward normal to $\partial\Omega$ at that point. For any compact $K\subset \R_-^N$, we have for $n$ large enough, that $\frac{\Omega}{r_n}\cap K=\emptyset$ as $r_n\rightarrow 0$. Since $supp\;v_n\subset \frac{\Omega}{r_n}$ and $v_n\rightarrow v$ a.e. in $\R^N$, it follows that $v=0$ a.e. on $K$. Therefore,  $supp\;v\subset \R_+^N$.
Hence, for $n$ large enough, we obtain that $supp\;v_n\subset \R_+^N$
and $v_n\rightharpoonup v$ in $D_{0}^{1,2}(\R_+^N)$.
We note that the functional $A_{\lambda_1,\lambda_2}$ is invariant under dilation, hence,
$$\|v_n\|^2=\int_{\R^N}|\nabla (r_{n}^{\frac{N-2}{2}}u_n(r_nx))|^2dx=\int_{\R^N} |\nabla u_n|^2 dx=\|u_n\|^2$$
and
$$\int_{\R^N}\frac{|v_n|^{2^*(s_i)}}{|x|^{s_i}}dx=
\int_{\R^N}r_{n}^{N-s_i}\frac{|u_n(r_nx)|^{2^*(s_i)}}{|x|^{s_i}}dx
=\int_{\R^N}\frac{|u_n|^{2^*(s_i)}}{|x|^{s_i}}dx.$$
Similarly, we have that
$\|w_n\|^2=\|(r_n)^{\frac{N-2}{2}}w_n(r_n x)\|^2.$
Notice that
$(r_n)^{\frac{N-2}{2}}w_n(r_n x)=v_n-v.$
When $v_n\rightharpoonup v$ in $D^{1,2}(\R^N)$, we have
$$\|w_n\|^2= \|v_n-v\|^2=\|v_n\|^2-\|v\|^2+o(1)
=\|u_n\|^2-\|v\|^2+o(1).$$
Recalling that $v_n\rightharpoonup v$ in $D^{1,2}(\R^N)$,
by  Brezis-Lieb type lemma (see \cite{BrezisLieb.1983} for s=0  and  \cite{GhoussoubKang.2004} for $s>0$)  and the invariance  property again,  we have
\begin{align*}
& A_{\lambda_1,\lambda_2}(w_n)\\
&=
\frac{1}{2}\int_{\R^N}|\nabla w_n|^2dx-\frac{\lambda_1}{2^*(s_1)}\int_{\R^N}\frac{|w_n|^{2^*(s_1)}}{|x|^{s_1}}dx
-\frac{\lambda_2}{2^*(s_2)}\int_{\R^N}\frac{|w_n|^{2^*(s_2)}}{|x|^{s_2}}dx\\
&=A_{\lambda_1,\lambda_2}\big(r_{n}^{\frac{N-2}{2}}w_n(r_nx)\big)\\
&=A_{\lambda_1,\lambda_2}(v_n-v)\\
&=A_{\lambda_1,\lambda_2}(v_n)-A_{\lambda_1,\lambda_2}(v)+o(1)\\
&=A_{\lambda_1,\lambda_2}(u_n)-A_{\lambda_1,\lambda_2}(v)+o(1)\\
&=c-A_{\lambda_1,\lambda_2}(v)+o(1).
\end{align*}
For any $h\in C_0^\infty(\R_+^N)$, let $h_n(x):=(r_n)^{\frac{2-N}{2}}h(\frac{x}{r_n})$, then we have that $h_n\in H_0^1(\Omega)$ for $n$ large enough due to the assumption that $r_n\rightarrow 0$. Thus
\begin{align*}
\langle A'_{\lambda_1,\lambda_2}(v), h\rangle=&\langle A'_{\lambda_1,\lambda_2}(v_n), h\rangle +o(1)\\
=&\langle A'_{\lambda_1,\lambda_2}(u_n), h_n\rangle +o(1)\\
=&o(1)\|h_n\|+o(1)\\
=&o(1)\|h\|+o(1),
\end{align*}
which implies that $A'_{\lambda_1,\lambda_2}(v)=0$.
For any $h\in H_0^1(\Omega)$, let $\tilde{h}_n(x):=r_{n}^{\frac{N-2}{2}}h(r_n x)$. Then for $n$ large enough, $supp\;\tilde{h}_n\subset \R_+^N$.
By  the Brezis-Lieb type lemma again,  we obtain that
\be\lab{2014-4-8-zhe1}
A'_{\lambda_1,\lambda_2}(v_n)-A'_{\lambda_1,\lambda_2}(v_n-v)-A'_{\lambda_1,\lambda_2}(v)\rightarrow 0\;\hbox{in}\;H^{-1}(\R^N).
\ee
Hence, for any $h\in H_0^1(\Omega)$,
\begin{align*}
&\langle A'_{\lambda_1,\lambda_2}(w_n), h\rangle=\langle A'_{\lambda_1,\lambda_2}(r_{n}^{\frac{N-2}{2}}w_n(r_nx)), \tilde{h}_n(x)\rangle\\
&=\langle A'_{\lambda_1,\lambda_2}(r_{n}^{\frac{N-2}{2}}w_n(r_nx)), \tilde{h}_n(x)\rangle+\langle A'_{\lambda_1,\lambda_2}(v(x)), \tilde{h}_n(x)\rangle\;\hbox{(since $A'_{\lambda_1,\lambda_2}(v)=0$)}\\
&=\langle A'_{\lambda_1,\lambda_2}(v_n-v), \tilde{h}_n(x)\rangle+\langle A'_{\lambda_1,\lambda_2}(v(x)), \tilde{h}_n(x)\rangle\\
&=\langle A'_{\lambda_1,\lambda_2}(v_n), \tilde{h}_n(x)\rangle+o(1)\|\tilde{h}_n\|\quad \hbox{(by (\ref{2014-4-8-zhe1}))}\\
&=\langle A'_{\lambda_1,\lambda_2}(u_n), h(x)\rangle+o(1)\|\tilde{h}_n\|\\
&=o(1)\|h\|\quad \hbox{(since $\|\tilde{h}_n\|\equiv\|h\|$)}.
\end{align*}
\ep

\bl\lab{2014-3-19-wl1}(See \cite[Lemma 3.5]{GhoussoubKang.2004})
If $u\in D^{1,2}(\R^N)$ and $h\in C_0^\infty(\R^N)$, then
$$\int_{\R^N}\frac{h^2|u|^{2^*(s)}}{|x|^s}dx\leq \mu_s(\R^N)^{-1}\Big(\int_{supp\;h}\frac{|u|^{2^*(s)}}{|x|^s}\Big)^{\frac{2^*(s)-2}{2^*(s)}}\int_{\R^N}|\nabla (hu)|^2dx,$$
where
\be\lab{2014-3-19-xwe1}
\mu_s(\R^N):=\inf\Big\{\int_{\R^N}|\nabla u|^2dx:     \quad u\in D_{0}^{1,2}(\R^N)\;\hbox{and}\;\int_{\R^N}\frac{|u|^{2^*(s)}}{|x|^s}<\infty\Big\}.
\ee
\el

\vskip 0.06in
\noindent{\bf Proof of Theorem  \ref{2014-3-19-compactness-th1}. }
Let $\{u_n\}\subset H_0^1(\Omega)$ be a bounded $(PS)_c$ sequence of $\Phi(u)$.
Up to a subsequence,     there is     an  $U^0\in H_0^1(\Omega)$ such that  $u_n\rightharpoonup U^0$ in $H_0^1(\Omega)$ and $\nabla u_n\rightarrow \nabla U^0$ a.e. on $\R^N$.  Evidently,  $\Phi'(U^0)=0$. Moreover, the sequence $u_n^1:=u_n-U^0$ satisfies
\be\lab{2014-3-19-we1}
\begin{cases}
\|u_n^1\|=\|u_n\|^2-\|U^0\|^2+o(1),\\
A'_{\lambda_1,\lambda_2}(u_n^1)\rightarrow 0\quad \hbox{in}\;H^{-1}(\Omega),\\
A_{\lambda_1,\lambda_2}(u_n^1)\rightarrow c-\Phi(U^0).
\end{cases}
\ee
If $u_n^1\rightarrow 0$ in $H_0^1(\Omega)$, we are done. If not, it is easy to see that
\be\lab{2014-3-19-we2}
\eta_0:=\liminf_{n\rightarrow \infty} \Big(\lambda_1\int_\Omega \frac{|u_n^1|^{2^*(s_1)}}{|x|^{s_1}}dx
+\lambda_2\int_\Omega \frac{|u_n^1|^{2^*(s_2)}}{|x|^{s_2}}dx\Big)
>0.
\ee
For the case of $\lambda_1>0$, we define an analogue of Levy's concentration function
$$Q_n(r):=\int_{B(0,r)}\Big(\lambda_1 \frac{|u_n^1|^{2^*(s_1)}}{|x|^{s_1}}
+\lambda_2\frac{|u_n^1|^{2^*(s_2)}}{|x|^{s_2}}\Big)dx.$$
Since $Q_n(0)=0$ and $Q_n(\infty)\geq \eta_0>0$, there exists a sequence $r_n^1>0$ such that for each $n$
\be\lab{2014-3-19-we3}
\delta=\int_{B(0,r_n^1)}\Big(\lambda_1 \frac{|u_n^1|^{2^*(s_1)}}{|x|^{s_1}}
+\lambda_2\frac{|u_n^1|^{2^*(s_2)}}{|x|^{s_2}}\Big)dx,
\ee
here we take $\delta$ so small   that
\be\lab{2014-3-19-we4}
\lambda_{1}^{\frac{2}{2^*(s_1)}}\mu_{s_1}(\R^N)^{-1}\delta^{\frac{2^*(s_1)-2}{2^*(s_1)}}
+\lambda_{2}^{\frac{2}{2^*(s_2)}}\mu_{s_2}(\R^N)^{-1}\delta^{\frac{2^*(s_2)-2}{2^*(s_2)}}
<\frac{1}{2},
\ee
where
$\mu_s(\R^N)$ is defined by (\ref{2014-3-19-xwe1}).   Define $v_n^1(x):=(r_n^1)^{\frac{N-2}{2}}u_n^1(r_n^1x)$. Since $\|v_n^1\|=\|u_n^1\|$ is bounded, we may assume $v_n^1\rightharpoonup U^1$ in $D^{1,2}(\R^N)$, $v_n^1\rightarrow U^1$ a.e. on $\R^N$ and
$$\delta=\int_{B(0,1)}\Big(\lambda_1 \frac{|v_n^1|^{2^*(s_1)}}{|x|^{s_1}}
+\lambda_2\frac{|v_n^1|^{2^*(s_2)}}{|x|^{s_2}}\Big)dx.$$
Next, we show that $U^1\not\equiv 0$.  Define $\Omega_n=\frac{1}{r_n^1}\Omega$ and let $f_n\in H_0^1(\Omega)$ be such that for any $h\in H_0^1(\Omega)$, we have $\displaystyle \langle A'_{\lambda_1,\lambda_2}(u_n^1), h\rangle=\int_\Omega \nabla f_n\cdot \nabla h$. Then $g_n(x):=(r_n^1)^{\frac{N-2}{2}}f_n(r_n^1x)$ satisfies
$\displaystyle \int_{\Omega_n}|\nabla g_n|^2 =\int_\Omega |\nabla f_n|^2$ and $\displaystyle \langle A'_{\lambda_1,\lambda_2}(v_n^1), h\rangle=\int_{\Omega_n}\nabla g_n\cdot \nabla h$ for any $h\in H_0^1(\Omega_m)$.   If $U_1\equiv 0$,  then for any $h\in C_0^\infty(\R^N)$  with  $supp\;h\subset B(0,1)$,   from Lemma \ref{2014-3-19-wl1} and the fact of (\ref{2014-3-19-we4}), we get that \begin{align*}
& \int_{B(0,1)}|\nabla (hv_n^1)|^2\\&=\int_{B(0,1)}\nabla v_n^1\cdot \nabla (h^2v_n^1)+o(1)\\
&=\lambda_1\int \frac{h^2|v_n^1|^{2^*(s_1)}}{|x|^{s_1}}+\lambda_2\int \frac{h^2|v_n^1|^{2^*(s_2)}}{|x|^{s_2}}+\int \nabla g_n\cdot \nabla (h^2v_n^1)+o(1)\\
&= \lambda_1\int \frac{h^2|v_n^1|^{2^*(s_1)}}{|x|^{s_1}}+\lambda_2\int \frac{h^2|v_n^1|^{2^*(s_2)}}{|x|^{s_2}}+\langle A'_{\lambda_1,\lambda_2}(v_n), h^2v_n^1\rangle+o(1)\\
& \leq \lambda_1 \mu_{s_1}(\R^N)^{-1}\Big(\int_{B(0,1)}\frac{|u_n^1|^{2^*(s_1)}}{|x|^{s_1}}\Big)^{\frac{2-2^*(s_1)}{2}}\int |\nabla (hv_n^1)|^2\\
&\quad +\lambda_2 \mu_{s_2}(\R^N)^{-1}\Big(\int_{B(0,1)}\frac{|u_n^1|^{2^*(s_2)}}{|x|^{s_2}}\Big)^{\frac{2-2^*(s_2)}{2}}\int |\nabla (hv_n^1)|^2+o(1)\\
& \leq \frac{1}{2}\int |\nabla (hv_n^1)|^2+o(1).
\end{align*}
Hence,  $\nabla v_n^1\rightarrow 0$ in $L_{loc}^{2}\big(B(0,1)\big)$ and $v_n^1\rightarrow 0$ in $L^{2^*(s_i)}_{loc}\big(B(0,1), |x|^{-s_i}dx\big)$, which contradicts the fact that
$$\int_{B(0,1)}\Big(\lambda_1 \frac{|v_n^1|^{2^*(s_1)}}{|x|^{s_1}}
+\lambda_2\frac{|v_n^1|^{2^*(s_2)}}{|x|^{s_2}}\Big)dx=\delta>0.$$
Thus we have proved that $U^1\not\equiv 0$.    Apply the similar argument for the case of $\lambda_1<0$ with a modified concentration function
 $$Q_n(r):=\int_{B(0,r)}\lambda_2\frac{|u_n^1|^{2^*(s_2)}}{|x|^{s_2}}dx.$$
 In this case ,we take $0<\delta<\big(\frac{\mu_{s_2}(\R^N)}{2}\big)^{\frac{N-s_2}{2-s_2}}$ small enough and  a sequence $r_n^1>0$ with $Q_{r_n^1}(x)=\delta$.
 We also define $v_n^1(x):=(r_n^1)^{\frac{N-2}{2}}u_n^1(r_n^1x)$ and assume that $v_n^1\rightharpoonup U^1$ in $D^{1,2}(\R^N)$, $v_n^1\rightarrow U^1$ a.e. on $\R^N$ and
$$\delta=\int_{B(0,1)}\lambda_2\frac{|v_n^1|^{2^*(s_2)}}{|x|^{s_2}}dx.$$
Next we will prove that $U^1\not\equiv 0$ for this case.

If $U_1\equiv 0$,  choose  any $h\in C_0^\infty(\R^N)$ such that $supp\;h\subset B(0,1)$  and invoke  Lemma \ref{2014-3-19-wl1} and the fact of $0<\delta<\big(\frac{\mu_{s_2}(\R^N)}{2}\big)^{\frac{N-s_2}{2-s_2}}$,  we have the following  estimates:
\begin{align*}
& \int_{B(0,1)}|\nabla (hv_n^1)|^2\\
&=\int_{B(0,1)}\nabla v_n^1\cdot \nabla (h^2v_n^1)+o(1)\\
&=\lambda_1\int \frac{h^2|v_n^1|^{2^*(s_1)}}{|x|^{s_1}}+\lambda_2\int \frac{h^2|v_n^1|^{2^*(s_2)}}{|x|^{s_2}}+\int \nabla g_n\cdot \nabla (h^2v_n^1)+o(1)\\
&\leq \lambda_2\int \frac{h^2|v_n^1|^{2^*(s_2)}}{|x|^{s_2}}+\int \nabla g_n\cdot \nabla (h^2v_n^1)+o(1)\\
& \leq \lambda_2 \mu_{s_2}(\R^N)^{-1}\Big(\int_{B(0,1)}\frac{|u_n^1|^{2^*(s_2)}}{|x|^{s_2}}\Big)^{\frac{2-2^*(s_2)}{2}}\int |\nabla (hv_n^1)|^2+o(1)\\
& \leq \frac{1}{2}\int |\nabla (hv_n^1)|^2+o(1).
\end{align*}
Hence,  $\nabla v_n^1\rightarrow 0$ in $L_{loc}^{2}\big(B(0,1)\big)$ and $v_n^1\rightarrow 0$ in $L^{2^*(s_2)}_{loc}\big(B(0,1), |x|^{-s_2}dx\big)$, which contradicts the fact that
$$\int_{B(0,1)}\Big(\lambda_2\frac{|v_n^1|^{2^*(s_2)}}{|x|^{s_2}}\Big)dx=\delta>0.$$
Thus $U^1\not\equiv 0$ is also true  for the case of $\lambda_1<0$.  In either case, we will prove that $r_n^1\rightarrow 0$. If not, since $\Omega$ is bounded, we may assume that $r_n^1\rightarrow r_\infty^1>0$, the fact that $u_n^1\rightharpoonup 0$ in $H_0^1(\Omega)$ means   that $v_n^1(x):=(r_n^1)^{\frac{N-2}{2}}u_n^1(r_n^1x)\rightharpoonup 0$ in $D_{0}^{1,2}(\R^N)$, which contradicts the fact $U^1\not\equiv 0$, and therefore $r_n^1\rightarrow 0$.

 Next, we prove that $supp\;U^1\subset \R_+^N$. Without loss of generality, assume that  $\partial R_+^N:=\{x_N=0\}$ is tangent to $\partial\Omega$ at $0$, and that $-e_N=(0,\cdots, -1)$ is the outward normal to $\partial\Omega$ at that point. For any compact $K\subset \R_-^N$, we have for $n$ large enough, that $\frac{\Omega}{r_n^1}\cap K=\emptyset$ as $r_n^1\rightarrow 0$. Since $supp\;v_n^1\subset \frac{\Omega}{r_n^1}$ and $v_n^1\rightarrow U^1$ a.e. in $\R^N$, it follows that $U^1=0$ a.e. on $K$, and therefore $supp\;U^1\subset \R_+^N$. By (\ref{2014-3-19-we1}) and Lemma \ref{2014-3-19-l2}, $A'_{\lambda_1,\lambda_2}(U^1)=0$ and $U^1$ is a weak solution of (\ref{2014-3-19-e5}). The sequence $u_n^2(x):=u_n^1(x)-(r_n^1)^{\frac{2-N}{2}}U^1(\frac{x}{r_n^1})$ also satisfies
\be\lab{2014-3-19-wse1}
\begin{cases}
\|u_n^2\|^2=\|u_n\|^2-\|U^0\|^2-\|U^1\|^2+o(1),\\
A_{\lambda_1,\lambda_2}(u_n^2)\rightarrow c-\Phi(U^0)-A_{\lambda_1, \lambda_2}(U^1),\\
A'_{\lambda_1,\lambda_2}(u_n^2)\rightarrow 0\;\hbox{in}\;H^{-1}(\Omega).
\end{cases}
\ee
Moreover,
$$A_{\lambda_1,\lambda_2}(U^1)\geq c_{\lambda_1,\lambda_2}>0.$$
By iterating the above procedure, we construct similarly sequences $U^j, (r_n^j)$ with the above properties  and $U^j$ is a solution of (\ref{2014-3-19-e5}).
It is easy to see that the iteration must   terminate after a finite number of  steps.\hfill$\Box$

\vskip 0.4in
Next, we will prove that $\displaystyle c_0:=\inf_{u\in \mathcal{N}}\Phi(u)<c_{\lambda_1,\lambda_2}$. Firstly, we recall  the following result.
\bl\lab{2014-3-24-l1}(\cite[Theorem 1.1]{LiLin.2012})
Suppose $\Omega$ is a bounded smooth domain in $\R^N$, $0\in \partial\Omega$ and the mean curvature $H(0)<0$. Then the equation
\be\lab{2014-3-24-e2}
\begin{cases}
\Delta u+\lambda \frac{u^{2^*(s_1)-1}}{|x|^{s_1}}+\frac{u^{2^*(s_2)-1}}{|x|^{s_2}}=0\quad &\hbox{in}\;\Omega,\\
u(x)>0\quad \hbox{in}\;\Omega, \quad u(x)=0\quad &\hbox{on}\;\partial\Omega,
\end{cases}
\ee
has a least-energy solution if  $N\geq 3, \lambda\in \R\;\hbox{and}\;0<s_2<s_1<2.$
\el
\br\lab{2014-3-24-r1}
Let $c_1$ be the least energy  corresponding to (\ref{2014-3-24-e2}). It has been proved that
$c_1<c_{\lambda,1}.$
We refer to \cite[Lemma 4.1]{LiLin.2012}.
\er

\bc\lab{2014-3-24-cro1}
Suppose $\Omega$ is a bounded smooth domain in $\R^N$, $0\in \partial\Omega$ and the mean curvature $H(0)<0$.  Assume $N\geq 3, 0<s_2<s_1<2, 0\leq s_3<2,\lambda_2,\lambda_3>0, 1<p<2^*(s_3)-1$. Furthermore, $\lambda_1>0$ or $\lambda_1<0$ with $p\geq 2^*(s_1)-1$,  then
\be\lab{2014-3-24-e3}
c_0:=\inf_{u\in \mathcal{N}}\Phi(u)<c_{\lambda_1,\lambda_2}.
\ee
\ec
\bp
It is easy to see that
 \be\lab{2014-3-25-e3}
\begin{cases}
\Delta u+\lambda_1 \frac{u^{2^*(s_1)-1}}{|x|^{s_1}}+\lambda_2\frac{u^{2^*(s_2)-1}}{|x|^{s_2}}=0\quad &\hbox{in}\;\Omega,\\
u(x)>0\quad \hbox{in}\;\Omega\quad u(x)=0\quad &\hbox{on}\;\partial\Omega
\end{cases}
\ee
 has a least-energy solution for $\lambda_2>0, N\geq 3, \lambda_1\in \R, 0<s_2<s_1<2 $ (thanks  to Lemma \ref{2014-3-24-l1}). For this case, we  denote the corresponding least energy by $\hat{c}_{\lambda_1,\lambda_2}$. Let $\lambda=\lambda_1 \lambda_{2}^{\frac{2-2^*(s_1)}{2^*(s_2)-2}}$, by (\ref{2014-3-19-e7}) and Remark \ref{2014-3-24-r1}, we have
\be\lab{2014-3-24-e5}
\hat{c}_{\lambda_1,\lambda_2}=\lambda_{2}^{\frac{-2}{2^*(s_2)-2}}c_1<\lambda_{2}^{\frac{-2}{2^*(s_2)-2}}c_{\lambda,1}
=c_{\lambda_1,\lambda_2}.
\ee
We note that  the assumptions required in Lemma \ref{Nehari-l1} are satisfied.  Let $w\in H_0^1(\Omega)$ be a least-energy solution to (\ref{2014-3-25-e3}). It is easy to see that
$J(w)=\max_{t>0}J(tw)=\hat{c}_{\lambda_1,\lambda_2},$
where
$$J(w):=\frac{1}{2}\int_{\Omega}|\nabla w|^2dx-\frac{\lambda_1}{2^*(s_1)}\int_{\Omega}\frac{|w|^{2^*(s_1)}}{|x|^{s_1}}dx
-\frac{\lambda_2}{2^*(s_2)}\int_{\Omega}\frac{|w|^{2^*(s_2)}}{|x|^{s_2}}dx.$$
Then for such a $w$, there exists some $t_w>0$ such that $t_ww\in \mathcal{N}$.  It follows that
\begin{align*}
c_0:=\inf_{u\in\mathcal{N}}\Phi(u)\leq \Phi(t_ww)
<J(t_ww)\leq J(w)=\hat{c}_{\lambda_1,\lambda_2}<c_{\lambda_1,\lambda_2}.
\end{align*}
\ep
However, for the case of $\lambda_3<0$, similar to the arguments of \cite[Lemma 3.1]{CeramiZhongZou.2014}, we can prove that $c_0\geq \hat{c}_{\lambda_1,\lambda_2}$. Nevertheless, the following lemma shows that $c_0<c_{\lambda_1,\lambda_2}$.
\bl\lab{2014-3-24-l2}
 Let $\Omega$ be a bounded smooth domain in $\R^N$ with $0\in \partial\Omega$. Suppose  that  $\partial\Omega$ is $C^2$ at $0$ and the mean curvature $H(0)<0$. We also assume that $N\geq 3, 0<s_2<s_1<2, 0\leq s_3<2, \lambda_2>0, \lambda_3<0, 1<p<2^*(s_3)-1$ and $\begin{cases} \lambda_1>0\\ p\leq 2^*(s_1)-1 \end{cases}$ or $\begin{cases}\lambda_1<0\\ p<2^*(s_2)-1\end{cases}$.
 Then we have
 $c_0:=\inf_{u\in\mathcal{N}}\Phi(u)<c_{\lambda_1,\lambda_2}$ if one of the following additional conditions  is satisfied:
 \begin{itemize}
 \item[(1)] $p<\frac{N-2s_3}{N-2}$.
 \item[(2)] $p\geq\frac{N-2s_3}{N-2}, |\lambda_3|$ is sufficiently small.
 \end{itemize}
\el

\bp
We prove this lemma by a modification of \cite[Lemma 2.2]{HsiaLinWadade.2010}.
Without loss of generality, we may assume that in a neighborhood of $0$, $\partial\Omega$ can be represented by $x_N=\varphi(x')$, where $x'=(x_1,\cdots, x_{N-1}), \varphi(0)=0, \nabla'\varphi(0)=0, \nabla'=(\partial_1, \cdots, \partial_{N-1})$ and the outer normal of $\partial\Omega$ at $0$ is $-e_N=(0,\cdots, 0, -1)$. Define $\phi(x)=(x', x_N-\varphi(x'))$ to ``flatten out" the boundary. We can choose a small $r_0>0$ and neighborhoods of $0$, $U$ and $\tilde{U}$, such that
$\displaystyle\phi(U)=B_{r_0}(0),\;\phi(U\cap \Omega)=B_{r_0}^{+}(0),
\phi(\tilde{U})=B_{\frac{r_0}{2}}(0)$ and $\phi(\tilde{U}\cap \Omega)=B_{\frac{r_0}{2}}^{+}(0)$. Here we adopt the notation:
$B_{r_0}^{+}(0)=B_{r_0}\cap \R_+^N\;\hbox{for any}\;r_0>0.$
Since $\partial\Omega\in C^2$, $\varphi$ can be expanded by
\be\lab{2014-3-1-e0}
\varphi(y')=\sum_{i=1}^{N-1}\alpha_i y_i^2+o(|y'|^2).
\ee
Then
$$H(0)=\frac{1}{N-1}\sum_{i=1}^{N-1}\alpha_i.$$
Suppose   that  $u\in H_0^1(\R_+^N)$ is a least-energy solution of (\ref{2014-3-19-e5}), i.e.,
\be\lab{2014-3-24-we1}
\begin{cases}
\Delta u+\lambda_1 \frac{u^{2^*(s_1)-1}}{|x|^{s_1}}+\lambda_2\frac{u^{2^*(s_2)-1}}{|x|^{s_2}}=0\quad &\hbox{in}\;\R_+^N,\\
u(x)>0\quad \hbox{in}\;\R_+^N,  \quad u(x)=0\quad &\hbox{on}\;\partial\R_+^N,
\end{cases}
\ee
and
$$A_{\lambda_1,\lambda_2}(u)=\frac{1}{2}a(u)-\frac{\lambda_1}{2^*(s_1)}b(u)-\frac{\lambda_2}{2^*(s_2)}c(u)=c_{\lambda_1,\lambda_2},$$
where $a(u), b(u), c(u)$ are defined by (\ref{2014-3-19-e1}). We also note that
\be\lab{2014-3-24-wbue1}
\max_{t>0}A_{\lambda_1,\lambda_2}(tu)=A_{\lambda_1,\lambda_2}(u)=c_{\lambda_1,\lambda_2}.
\ee
Let $\varepsilon>0$, we define
$$v_\varepsilon(x):=\varepsilon^{-\frac{N-2}{2}}u(\frac{\phi(x)}{\varepsilon})\;\hbox{for}\;x\in \Omega\cap U.$$
Let $\eta\in C_0^\infty(U)$ be a positive cut-off function with $\eta\equiv 1$ in $\tilde{U}$ and consider $\hat{v}_\varepsilon:=\eta v_\varepsilon$ in $\Omega$, then for $t\geq 0$, if $\lambda_1>0$, we have
\begin{align*}
\Phi(t\hat{v}_\varepsilon)=&\frac{t^2}{2}\int_\Omega |\nabla \hat{v}_\varepsilon|^2 dx -\lambda_1\frac{ t^{2^*(s_1)}}{2^*(s_1)}\int_\Omega \frac{\hat{v}_{\varepsilon}^{2^*(s_1)}}{|x|^{s_1}}dx\\
&-\lambda_2\frac{t^{2^*(s_2)}}{2^*(s_2)}
\int_\Omega\frac{\hat{v}_{\varepsilon}^{2^*(s_2)}}{|x|^{s_2}}dx-\lambda_3\frac{ t^{p+1}}{p+1}\int_\Omega \frac{\hat{v}_{\varepsilon}^{p+1}}{|x|^{s_3}}dx\\
\leq&\frac{t^2}{2}\int_\Omega |\nabla \hat{v}_\varepsilon|^2 dx -\lambda_1\frac{t^{2^*(s_1)}}{2^*(s_1)}
\int_{\Omega\cap \tilde{U}}\frac{v_{\varepsilon}^{2^*(s_1)}}{|x|^{s_1}}dx\\
&-\lambda_2\frac{t^{2^*(s_2)}}{2^*(s_2)}
\int_{\Omega\cap \tilde{U}}\frac{v_{\varepsilon}^{2^*(s_2)}}{|x|^{s_2}}dx-\lambda_3\frac{ t^{p+1}}{p+1}\int_\Omega \frac{\hat{v}_{\varepsilon}^{p+1}}{|x|^{s_3}}dx.
\end{align*}
If $\lambda_1<0$, we have
\begin{align*}
\Phi(t\hat{v}_\varepsilon)\leq&\frac{t^2}{2}\int_\Omega |\nabla \hat{v}_\varepsilon|^2 dx -\lambda_1\frac{t^{2^*(s_1)}}{2^*(s_1)}
\int_{\Omega\cap U}\frac{v_{\varepsilon}^{2^*(s_1)}}{|x|^{s_1}}dx\\
&-\lambda_2\frac{t^{2^*(s_2)}}{2^*(s_2)}
\int_{\Omega\cap \tilde{U}}\frac{v_{\varepsilon}^{2^*(s_2)}}{|x|^{s_2}}dx-\lambda_3\frac{ t^{p+1}}{p+1}\int_\Omega \frac{\hat{v}_{\varepsilon}^{p+1}}{|x|^{s_3}}dx.
\end{align*}
By the change of the variable $y=\frac{\phi(x)}{\varepsilon}$, we have
\begin{align*}
\int_\Omega |\nabla \hat{v}_\varepsilon|^2dx=&\int_{\Omega\cap U}\eta^2|\nabla v_\varepsilon|^2dx-\int_{\Omega\cap U}\eta(\Delta \eta)v_\varepsilon^2 dx\\
\leq&\int_{\R_+^N}|\nabla u(y)|^2 dy-2\int_{B_{\frac{r_0}{\varepsilon}}^{+}}\eta\big(\phi^{-1}(\varepsilon y)\big)^2\partial_Nu(y)\nabla'u(y)\cdot (\nabla'\varphi)(\varepsilon y')dy\\
&+\int_{B_{\frac{r_0}{\varepsilon}}^{+}}\eta\big(\phi^{-1}(\varepsilon y)\big)^2\big(\partial_Nu(y)\big)^2|(\nabla'u)(\varepsilon y')|^2dy\\
&-\varepsilon^2\int_{B_{\frac{r_0}{\varepsilon}}^{+}}\eta\big(\phi^{-1}(\varepsilon y)\big)(\Delta \eta)\big(\phi^{-1}(\varepsilon y)\big)u(y)^2 dy.
\end{align*}
Note that, by using $|\nabla' \phi(y')|=O(|y'|)$ and the decay estimate of $|\nabla u|$ in (\ref{2014-3-19-e9}), we see that
$$\int_{B_{\frac{r_0}{\varepsilon}}^{+}}\eta\big(\phi^{-1}(\varepsilon y)\big)^2\big(\partial_Nu(y)\big)^2|(\nabla'u)(\varepsilon y')|^2dy\leq C\varepsilon^2\int_{\R^N}\big(1+|y|\big)^{-2N}|y|^2dy=O(\varepsilon^2).$$
Hence,
$$\int_\Omega |\nabla \hat{v}_\varepsilon|^2dx=a(u)-2\int_{B_{\frac{r_0}{\varepsilon}}^{+}}\eta\big(\phi^{-1}(\varepsilon y)\big)^2\partial_Nu(y)\nabla'u(y)\cdot (\nabla'\varphi)(\varepsilon y')dy+O(\varepsilon^2).$$
Using integration by parts and the formulas (\ref{2014-3-19-e9}), (\ref{2014-3-1-e0}),  we see that
\begin{align*}
I:=&-2\int_{B_{\frac{r_0}{\varepsilon}}^{+}}\eta\big(\phi^{-1}(\varepsilon y)\big)^2\partial_Nu(y)\nabla'u(y)\cdot (\nabla'\varphi)(\varepsilon y')dy\\
=&-\frac{2}{\varepsilon}\int_{B_{\frac{r_0}{\varepsilon}}^{+}}\eta\big(\phi^{-1}(\varepsilon y)\big)^2\partial_Nu(y)\nabla'u(y)\cdot \nabla'[\varphi(\varepsilon y')]dy\\
=&-\frac{2}{\varepsilon}\int_{B_{\frac{r_0}{\varepsilon}}^{+}\cap \partial\R_+^N}\eta\big(\phi^{-1}(\varepsilon y)\big)^2 \partial_N u(y)\nabla'u(y)\varphi(\varepsilon y')dS_y\\
&+\frac{4}{\varepsilon}\int_{B_{\frac{r_0}{\varepsilon}}^{+}}\eta\big(\phi^{-1}(\varepsilon y)\big)\nabla'[\eta\big(\phi^{-1}(\varepsilon y)\big)]\partial_Nu(y)\nabla'u(y)\cdot \varphi(\varepsilon y')dy\\
&+\frac{2}{\varepsilon}\int_{B_{\frac{r_0}{\varepsilon}}^{+}}\eta\big(\phi^{-1}(\varepsilon y)\big)^2\nabla'\partial_Nu(y)\nabla'u(y)\cdot \varphi(\varepsilon y')dy\\
&+\frac{2}{\varepsilon}\int_{B_{\frac{r_0}{\varepsilon}}^{+}}\eta\big(\phi^{-1}(\varepsilon y)\big)^2\partial_Nu(y)\sum_{i=1}^{N-1}\partial_{ii}u(y)\varphi(\varepsilon y')dy\\
=&\frac{2}{\varepsilon}\int_{B_{\frac{r_0}{\varepsilon}}^{+}}\eta\big(\phi^{-1}(\varepsilon y)\big)^2\partial_Nu(y)\sum_{i=1}^{N-1}\partial_{ii}u(y)\varphi(\varepsilon y')dy+O(\varepsilon^2).
\end{align*}
Applying (\ref{2014-3-24-we1}) and integration by parts, we obtain that
\begin{align*}
I':=&\frac{2}{\varepsilon}\int_{B_{\frac{r_0}{\varepsilon}}^{+}}\eta\big(\phi^{-1}(\varepsilon y)\big)^2\partial_Nu(y)\sum_{i=1}^{N-1}\partial_{ii}u(y)\varphi(\varepsilon y')dy\\
=&\frac{2}{\varepsilon}\int_{B_{\frac{r_0}{\varepsilon}}^{+}}\eta\big(\phi^{-1}(\varepsilon y)\big)^2\partial_Nu(y)[\Delta u(y)-\partial_{NN}u(y)]\varphi(\varepsilon y')dy\\
=&-\frac{2}{\varepsilon}\int_{B_{\frac{r_0}{\varepsilon}}^{+}}\eta\big(\phi^{-1}(\varepsilon y)\big)^2\partial_Nu(y)[\lambda_1\frac{u^{2^*(s_1)-1}}{|y|^{s_1}}+\lambda_2\frac{u^{2^*(s_2)-1}}{|y|^{s_2}}]\varphi(\varepsilon y')dy\\
&-\frac{1}{\varepsilon}\int_{B_{\frac{r_0}{\varepsilon}}^{+}}\eta\big(\phi^{-1}(\varepsilon y)\big)^2\partial_N[\big(\partial_N u(y)\big)^2]\varphi(\varepsilon y')dy\\
=&-\frac{2}{\varepsilon}\lambda_1\frac{1}{2^*(s_1)}\int_{B_{\frac{r_0}{\varepsilon}}^{+}}\eta\big(\phi^{-1}(\varepsilon y)\big)^2\frac{\partial_N[u(y)^{2^*(s_1)}]}{|y|^{s_1}}\varphi(\varepsilon y')dy\\
&-\frac{2}{\varepsilon}\lambda_2\frac{1}{2^*(s_2)}\int_{B_{\frac{r_0}{\varepsilon}}^{+}}\eta\big(\phi^{-1}(\varepsilon y)\big)^2\frac{\partial_N[u(y)^{2^*(s_2)}]}{|y|^{s_2}}\varphi(\varepsilon y')dy\\
&+\frac{1}{\varepsilon}\int_{B_{\frac{r_0}{\varepsilon}}^{+}\cap \partial \R_+^N}\eta\big(\phi^{-1}(\varepsilon y)\big)^2\big(\partial_N u(y)\big)^2\varphi(\varepsilon y')dS_y+O(\varepsilon^2)\\
=&-\frac{2}{\varepsilon}\lambda_1\frac{s_1}{2^*(s_1)}\int_{B_{\frac{r_0}{\varepsilon}}^{+}}\eta\big(\phi^{-1}(\varepsilon y)\big)^2\frac{u(y)^{2^*(s_1)}y_N}{|y|^{s_1+2}}\varphi(\varepsilon y')dy\\
&-\frac{2}{\varepsilon}\lambda_2\frac{s_2}{2^*(s_2)}\int_{B_{\frac{r_0}{\varepsilon}}^{+}}\eta\big(\phi^{-1}(\varepsilon y)\big)^2\frac{u(y)^{2^*(s_2)}y_N}{|y|^{s_2+2}}\varphi(\varepsilon y')dy\\
&+\frac{1}{\varepsilon}\int_{B_{\frac{r_0}{\varepsilon}}^{+}\cap \partial \R_+^N}\eta\big(\phi^{-1}(\varepsilon y)\big)^2\big(\partial_N u(y)\big)^2\varphi(\varepsilon y')dS_y+O(\varepsilon^2)\\
=:&J_1+J_2+J_3+O(\varepsilon^2).
\end{align*}
Among them
\begin{align*}
J_1:=&-\frac{2}{\varepsilon}\lambda_1\frac{s_1}{2^*(s_1)}\int_{B_{\frac{r_0}{\varepsilon}}^{+}}\eta\big(\phi^{-1}(\varepsilon y)\big)^2\frac{u(y)^{2^*(s_1)}y_N}{|y|^{s_1+2}}\varphi(\varepsilon y')dy\\
=&-\frac{2}{\varepsilon}\lambda_1\frac{s_1}{2^*(s_1)}\int_{B_{\frac{r_0}{\varepsilon}}^{+}\backslash B_{\frac{r_0/2}{\varepsilon}}^{+}}\eta\big(\phi^{-1}(\varepsilon y)\big)^2\frac{u(y)^{2^*(s_1)}y_N}{|y|^{s_1+2}}\varphi(\varepsilon y')dy\\
&-\frac{2}{\varepsilon}\lambda_1\frac{s_1}{2^*(s_1)}\int_{B_{\frac{r_0/2}{\varepsilon}}^{+}}\frac{u(y)^{2^*(s_1)}y_N}{|y|^{s_1+2}}\varphi(\varepsilon y')dy\\
=:&J_{1,1}+J_{1,2},
\end{align*}
and
$$|J_{1,1}|\leq C\varepsilon \int_{\{r_0/2\leq |\varepsilon y|<r_0\}}|y|^{2^*(s_1)(1-N)+1-s_1}dy=O(\varepsilon^{\frac{N(N-s_1)}{N-2}}).$$
Notice that
\be\lab{2014-3-24-we2}
\varepsilon \int_{\R_+^N\backslash B_{\frac{r_0/2}{\varepsilon}}^{+}}u(y)^{2^*(s_1)}|y|^{1-s_1}dy=O(\varepsilon^{\frac{N(N-s_1)}{N-2}}).
\ee
By (\ref{2014-3-1-e0}),  (\ref{2014-3-24-we2}) and using the fact of $u(y', y_N)=u(|y'|, y_N)$, we obtain
\begin{align*}
& J_{1,2}\\
&=-2\varepsilon\lambda_1\frac{s_1}{2^*(s_1)}\sum_{i=1}^{N-1}\alpha_i\int_{\R_+^N}\frac{u(y)^{2^*(s_1)}y_N}{|y|^{s_1+2}}y_i^2(1+o(1))dy+
O(\varepsilon^{\frac{N(N-s_1)}{N-2}})\\
&=-\frac{2\lambda_1s_1\varepsilon}{2^*(s_1)(N-1)}\int_{\R_+^N}\frac{u(y)^{2^*(s_1)}y_N}{|y|^{s_1+2}}|y'|^2dy\Big(\sum_{i=1}^{N-1}\alpha_i\Big)(1+o(1))
+
O(\varepsilon^{\frac{N(N-s_1)}{N-2}}).
\end{align*}
Thus,
$$
J_1=-\frac{2\lambda_1s_1}{2^*(s_1)}K_1H(0)\big(1+o(1)\big)\varepsilon+O(\varepsilon^2),
$$
where
$$K_1:=\int_{\R_+^N}\frac{u(y)^{2^*(s_1)}y_N}{|y|^{s_1+2}}|y'|^2dy.$$
Similarly, we can prove that
$$
J_2=-\frac{2\lambda_2s_2}{2^*(s_2)}K_2H(0)\big(1+o(1)\big)\varepsilon+O(\varepsilon^2),
$$
where
$$K_2:=\int_{\R_+^N}\frac{u(y)^{2^*(s_2)}y_N}{|y|^{s_2+2}}|y'|^2dy.$$
Next, we see that
\begin{align*}
J_3=&\frac{1}{\varepsilon}\int_{B_{\frac{r_0}{\varepsilon}}^{+}\cap \partial \R_+^N}\eta\big(\phi^{-1}(\varepsilon y)\big)^2\big(\partial_N u(y)\big)^2\varphi(\varepsilon y')dS_y\\
=&\frac{1}{\varepsilon}\int_{\big(B_{\frac{r_0}{\varepsilon}}^{+}\backslash B_{\frac{r_0/2}{\varepsilon}}^{+}\big)\cap \partial \R_+^N}\eta\big(\phi^{-1}(\varepsilon y)\big)^2\big(\partial_N u(y)\big)^2\varphi(\varepsilon y')dS_y\\
&+\frac{1}{\varepsilon}\int_{B_{\frac{r_0/2}{\varepsilon}}^{+}\cap \partial \R_+^N}\eta\big(\phi^{-1}(\varepsilon y)\big)^2\big(\partial_N u(y)\big)^2\varphi(\varepsilon y')dS_y\\
=:&J_{3,1}+J_{3,2},
\end{align*}
By the  mean value theorem for integrals, we have
\begin{align*}
|J_{3,1}|\leq& \frac{C(r_0)}{\varepsilon}\big(\frac{1}{\varepsilon}\big)^{-2N}\big(\frac{1}{\varepsilon}\big)^{N-1}\\
=&O(\varepsilon^{2N-1-(N-1)})=O(\varepsilon^N).
\end{align*}
Using the symmetry, by the polar coordinates transformation,
we also obtain that
\be\lab{2014-3-24-we3}
\varepsilon\int_{\{|\varepsilon y'|>\frac{r_0}{2}\}}|(\partial_N u)(y',0)|^2 |y'|^2dy'=O(\varepsilon^N).
\ee
Thus, by (\ref{2014-3-1-e0}), (\ref{2014-3-24-we3}) and using the fact of $u(y', y_N)=u(|y'|, y_N)$, we obtain
\begin{align*}
J_{3,2}=&\varepsilon \sum_{i=1}^{N-1}\alpha_i \int_{\R^{N-1}}\big((\partial_Nu)(y',0)\big)^2y_i^2dy'\big(1+o(1)\big)+O(\varepsilon^N)\\
=&\frac{\varepsilon}{N-1}\int_{\R^{N-1}}\big((\partial_Nu)(y',0)\big)^2|y'|^2dy'\sum_{i=1}^{N-1}\alpha_i+O(\varepsilon^2)\\
=&K_3H(0)\big(1+o(1)\big)\varepsilon+O(\varepsilon^2),
\end{align*}
where
$$K_3:=\int_{\R^{N-1}}\big((\partial_Nu)(y',0)\big)^2|y'|^2dy'>0.$$
Hence,
$$I'=\big(K_3-\frac{2\lambda_1s_1}{2^*(s_1)}K_1-\frac{2\lambda_2s_2}{2^*(s_2)}K_2\big)H(0)\big(1+o(1)\big)\varepsilon+O(\varepsilon^2),$$
which implies that
$$I=\big(K_3-\frac{2\lambda_1s_1}{2^*(s_1)}K_1-\frac{2\lambda_2s_2}{2^*(s_2)}K_2\big)H(0)\big(1+o(1)\big)\varepsilon+O(\varepsilon^2)$$
and that
$$\int_\Omega |\nabla \hat{v}_\varepsilon|^2dx=
a(u)+\big(K_3-\frac{2\lambda_1s_1}{2^*(s_1)}K_1-\frac{2\lambda_2s_2}{2^*(s_2)}K_2\big)H(0)\big(1+o(1)\big)\varepsilon+O(\varepsilon^2).$$
Furthermore, the integrals $\int_{\Omega \cap \tilde{U}}\frac{v_{\varepsilon}^{2^*(s_2)}}{|x|^{s_2}}dx,
\int_{\Omega \cap \tilde{U}}\frac{v_{\varepsilon}^{2^*(s_1)}}{|x|^{s_1}}dx$ and $
\int_{\Omega \cap U}\frac{v_{\varepsilon}^{2^*(s_1)}}{|x|^{s_1}}dx$ can be estimated  by  the same argument  as  that in  \cite[Lemma 2.2]{HsiaLinWadade.2010} to  obtain that
$$\int_{\Omega\cap \tilde{U}}\frac{v_{\varepsilon}^{2^*(s_2)}}{|x|^{s_2}}dx=c(u)-s_2K_2H(0)\big(1+o(1)\big)\varepsilon+O(\varepsilon^2)$$
and
$$\int_\Omega \frac{(\hat{v}_\varepsilon)^{2^*(s_1)}}{|x|^{s_1}}dx=b(u)-s_1K_1H(0)\big(1+o(1)\big)\varepsilon+O(\varepsilon^2).$$
By \cite[Lemma 2.4]{CeramiZhongZou.2014}, we also obtain that
$$\int_\Omega \frac{(\hat{v}_\varepsilon)^{p+1}}{|x|^{s_3}}dx=\varepsilon^{s_0-s_3}\int_{\R_+^N}\frac{u^{p+1}}{|y|^{s_3}}dy\big(1+o(1)\big)+O(\varepsilon^{\frac{N(p+1)}{2}}),$$
where
$s_0:=\frac{N+2}{2}-\frac{N-2}{2}p\in (s_3, 2)$ when $1<p<2^*(s_3)-1$.
Thus, we have
\begin{align*}
\Phi(t\hat{v}_\varepsilon)\leq&\frac{t^2}{2}\big[a(u)+\big(K_3-\frac{2\lambda_1s_1}{2^*(s_1)}K_1-\frac{2\lambda_2s_2}{2^*(s_2)}K_2\big)H(0)\big(1+o(1)\big)\varepsilon+O(\varepsilon^2)\big]\\
&-\lambda_1\frac{t^{2^*(s_1)}}{2^*(s_1)}\big[b(u)-s_1K_1H(0)\big(1+o(1)\big)\varepsilon+O(\varepsilon^2)\big]\\
&-\lambda_2\frac{t^{2^*(s_2)}}{2^*(s_2)}\big[c(u)-s_2K_2H(0)\big(1+o(1)\big)\varepsilon+O(\varepsilon^2)\big]\\
&-\lambda_3\frac{t^{p+1}}{p+1}\big[\int_{\R_+^N}\frac{u^{p+1}}{|y|^{s_3}}dy\big(1+o(1)\big)\varepsilon^{s_0-s_3}+O(\varepsilon^{\frac{N(p+1)}{2}})\big].
\end{align*}
Then, it is easy to see that there exists some $T$ large enough and $\varepsilon_0$ sufficiently small  such that  $\Phi(T\hat{v}_\varepsilon)<0$ for all $\varepsilon<\varepsilon_0$.
By Lemma \ref{Nehari-l1}  and Remark \ref{2014-3-18-r1},
  there exists a  unique $t_{\hat{v}_\varepsilon}>0$ such that $t_{\hat{v}_\varepsilon}\hat{v}_\varepsilon\in \mathcal{N}$ and
$ \max_{t>0}\Phi(t\hat{v}_\varepsilon)=\Phi(t_{\hat{v}_\varepsilon}\hat{v}_\varepsilon)\geq c_0>0.$
Hence, we obtain that
\be\lab{2014-3-24-wbue2}
t_{\hat{v}_\varepsilon}<T\;\hbox{for all}\;\varepsilon<\varepsilon_0.
\ee
On the other hand, when $t\in (0,T)$, we have
\begin{align*}
\Phi(t\hat{v}_\varepsilon)=&A_{\lambda_1,\lambda_2}(tu)+\big[\frac{t^2}{2}\big(K_3-\frac{2\lambda_1s_1}{2^*(s_1)}K_1-\frac{2\lambda_2s_2}{2^*(s_2)}K_2\big)\\
&+\lambda_1s_1K_1\frac{t^{2^*(s_1)}}{2^*(s_1)}+\lambda_2s_2K_2\frac{t^{2^*(s_2)}}{2^*(s_2)}\big]H(0)\big(1+o(1)\big)\varepsilon\\
&-\lambda_3\frac{t^{p+1}}{p+1}\int_{\R_+^N}\frac{u^{p+1}}{|y|^{s_3}}dy\big(1+o(1)\big)\varepsilon^{s_0-s_3}+O(\varepsilon^2)\\
=:&A_{\lambda_1,\lambda_2}(tu)+g_1(t)H(0)\big(1+o(1)\big)\varepsilon\\
&-\lambda_3\frac{t^{p+1}}{p+1}\int_{\R_+^N}\frac{u^{p+1}}{|y|^{s_3}}dy\big(1+o(1)\big)\varepsilon^{s_0-s_3}+O(\varepsilon^2).
\end{align*}
Notice that
$g_1(1)=\frac{1}{2}K_3>0,$
there exists $\delta_0>0,\varepsilon_1<\varepsilon_0$ such that
$$g_1(t)\geq \frac{1}{4}K_3>0\;\hbox{for all}\;(t,\varepsilon)\in [1-\delta_0,1+\delta_0]\times (0, \varepsilon_1).$$
Define
$$M:=\max_{t\in (0, 1-\delta_0]\cup [1+\delta_0, T]}A_{\lambda_1,\lambda_2}(tu),$$
we see  that
$M<c_{\lambda_1,\lambda_2}$ by  (\ref{2014-3-24-wbue1}).
Then, by the continuity, we obtain that for $\varepsilon<\varepsilon_2$ small enough,
\be\lab{2014-3-24-wbue3}
\max_{t\in (0, 1-\delta_0]\cup [1+\delta_0, T]}\Phi(t\hat{v}_\varepsilon)\leq M+O(\varepsilon^\sigma)<c_{\lambda_1,\lambda_2},
\ee
where $\sigma:=\min\{s_0-s_3, 1\}>0$.
On the other hand, recalling that $H(0)<0, g_1(t)>\frac{1}{4}K_3$ for all $1-\delta_0\leq t\leq 1+\delta_0$.
When $s_0-s_3>1$, that is, $p<\frac{N-2s_3}{N-2}$, then it is easy to see that there exists $\varepsilon_3<\varepsilon_2$ such that
$$g_1(t)H(0)\big(1+o(1)\big)\varepsilon-\lambda_3\frac{t^{p+1}}{p+1}
\int_{\R_+^N}\frac{u^{p+1}}{|y|^{s_3}}dy\big(1+o(1)\big)
\varepsilon^{s_0-s_3}+O(\varepsilon^2)<0$$
for all $(t,\varepsilon)\in [1-\delta_0, 1+\delta_0]\times (0, \varepsilon_3)$.
It follows that
\be\lab{2014-3-24-wbue4}
\max_{1-\delta_0\leq t\leq 1+\delta_0}\Phi(t\hat{v}_\varepsilon)<\max_{1-\delta_0\leq t\leq 1+\delta_0}A_{\lambda_1, \lambda_2}(tu)=A_{\lambda_1,\lambda_2}(u)=c_{\lambda_1,\lambda_2}.
\ee
Combine  with (\ref{2014-3-24-wbue2}), (\ref{2014-3-24-wbue3}) and (\ref{2014-3-24-wbue4}), we obtain that
$$c_0\leq \max_{t>0}\Phi(t\hat{v}_\varepsilon)<c_{\lambda_1,\lambda_2}$$
for $\varepsilon$ small enough.  Similar to Corollary \ref{2014-3-24-cro1}, let $w\in H_0^1(\Omega)$ be a least-energy solution to (\ref{2014-3-25-e3}). Then $$J(w)=\max_{t>0}J(tw)=\hat{c}_{\lambda_1,\lambda_2}<c_{\lambda_1,\lambda_2},$$
where
$$J(w):=\frac{1}{2}\int_{\Omega}|\nabla w|^2dx-\frac{\lambda_1}{2^*(s_1)}\int_{\Omega}\frac{|w|^{2^*(s_1)}}{|x|^{s_1}}dx
-\frac{\lambda_2}{2^*(s_2)}\int_{\Omega}\frac{|w|^{2^*(s_2)}}{|x|^{s_2}}dx.$$
Since $2^*(s_2)>2^*(s_1)$ and $\lambda_2>0$, it is easy to see that there exists some $T>0$ such that $J(Tw)<0$. When $|\lambda_3|$ is small enough, we have $\Phi(Tw)<0$.
By Lemma \ref{Nehari-l1} again, there exists some $0<t_w<T$ such that $t_ww\in \mathcal{N}$ and
$$\max_{t>0}\Phi(tw)=\Phi(t_ww)\geq c_0>0.$$
On the other hand, since $\hat{c}_{\lambda_1,\lambda_2}<c_{\lambda_1,\lambda_2}$, when $|\lambda_3|$ is small enough for the case of $p\geq \frac{N-2s_3}{N-2}$, we have
\begin{align*}
\Phi(t_ww)=&J(t_ww)-\lambda_3\frac{t_{w}^{p+1}}{p+1}\int_\Omega \frac{w^{p+1}}{|x|^{s_3}}dx\\
\leq&J(w)+|\lambda_3||T|^{p+1}\int_\Omega \frac{w^{p+1}}{|x|^{s_3}}dx\\
=&\hat{c}_{\lambda_1,\lambda_2}+|\lambda_3||T|^{p+1}\int_\Omega \frac{w^{p+1}}{|x|^{s_3}}dx\\
<&c_{\lambda_1,\lambda_2}.
\end{align*}
\ep


\vskip 0.13in
\noindent{\bf Proof of Theorem  \ref{main-th}.   }
Under the assumptions of Theorem \ref{main-th}, the Nehari manifold is well defined due to Lemma  \ref{Nehari-l1}.
By Ekeland's  variational principle,
let $\{u_n\}\subset \mathcal{N}$ be a minimizing sequence such that
$\Phi(u_n)\rightarrow c_0:=\inf_{u\in\mathcal{N}}\Phi(u)$
and $\Phi'|_{\mathcal{N}}(u_n)\rightarrow 0.$
By Lemma \ref{2014-3-25-l1} and Lemma  \ref{l3}, we obtain that  $\{u_n\}$ is  a bounded $(PS)_{c_0}$ sequence of  $\Phi$.  By Corollary  \ref{2014-3-24-cro1} or Lmema  \ref{2014-3-24-l2}, we can also obtain that
$c_0<c_{\lambda_1,\lambda_2}. $
Hence, by Corollary \ref{2014-3-19-cro1}, $\Phi(u)$ satisfies $(PS)_{c_0}$ condition. That is, up to a subsequence,    $u_n\rightarrow u_0$ strongly in $H_0^1(\Omega)$ for some $u_0\in H_0^1(\Omega)$ and $\Phi'(u_0)=0, \Phi(u_0)=c_0$. Thereby, the existence of  the ground state solution is established.
We also note that $\Phi(u)$ is even, which implies that $|u_0|\in \mathcal{N}$ and $\Phi(|u_0|)=\Phi(u_0)=c_0$. Hence, without loss of generality, we may assume that $u_0\geq 0$. Finally, apply the similar arguments as \cite[Proof of Lemma 2.6(i)]{LinWadadeothers.2012}, we can obtain the similar regularity property for a nonnegative solution. By the maximum principle, we claim  that $u>0$.\hfill$\Box$.


\s{Proof of Theorem  \ref{main-th2}}
\renewcommand{\theequation}{4.\arabic{equation}}
\renewcommand{\theremark}{4.\arabic{remark}}
\renewcommand{\thedefinition}{4.\arabic{definition}}
In this section, we always assume that $\lambda_3<0, p\leq 2^*(s_3)-1$.   Denote $\lambda_3$    by    $\lambda$ for  the simplicity. To obtain a positive solution, we consider the following modified functional
$$I_\lambda(u):=\frac{1}{2}\int_\Omega |\nabla u|^2dx
-\frac{\lambda_1}{2^*(s_1)}\int_\Omega \frac{u_{+}^{2^*(s_1)}}{|x|^{s_1}}dx$$
$$\quad\quad\quad\quad \quad\quad     -\frac{\lambda_2}{2^*(s_2)}\int_\Omega \frac{u_{+}^{2^*(s_2)}}{|x|^{s_2}}dx
-\frac{\lambda}{p+1}\int_\Omega \frac{u_{+}^{p+1}}{|x|^{s_3}}dx.$$
 By (\ref{2014-3-24-e5}), we have
$\hat{c}_{\lambda_1,\lambda_2}<c_{\lambda_1,\lambda_2}$ and it is easy to see that $\hat{c}_{\lambda_1,\lambda_2}$ is the ground state value of $I_0$ (i.e., $\lambda=0$) which can be obtained  by some $0<u_0\in H_0^1(\Omega)$. That is,
$$\Delta u_0+\lambda_1\frac{u_{0}^{2^*(s_1)-1}}{|x|^{s_1}}+\lambda_2\frac{u_{0}^{2^*(s_2)-1}}{|x|^{s_2}}=0$$
in $\Omega$ and
$I_0(u_0)=\hat{c}_{\lambda_1,\lambda_2}.$
Next, for the  convenience, we denote $\hat{c}_{\lambda_1,\lambda_2}$ by $c_0$.
It is easy to prove that $u_0$ is a mountain pass type solution. Here we  list the conditions which are fulfilled by $I_0$ without proof  (see \cite[Section 5]{CeramiZhongZou.2014}):
\begin{itemize}
\item[(M1)] there exists $c, r>0$ such that if $\|u\|=r$, then $I_0(u)\geq c$ and there exists $v_0\in H_0^1(\Omega)$ such that $\|v_0\|>r$ and $I_0(v_0)<0$;
\item[(M2)] there exists a critical point $u_0\in H_0^1(\Omega)$ of $I_0$ such that
    $$I_0(u_0)=c_0:=\min_{\gamma\in \Gamma}\max_{t\in [0,1]}I_0(\gamma(t)),$$
    where $\Gamma=\{\gamma\in C([0,1], H_0^1(\Omega))|\gamma(0)=0, \gamma(1)=v_0\}$;
\item[(M3)] it holds that
$c_0:=\inf\{I_0(u)|u\in H_0^1(\Omega)\backslash \{0\}\};$
\item[(M4)] the set $\mathcal{S}:=\{u\in H_0^1(\Omega)|I'_0(u)=0, I_0(u)=c_0\}$ is compact in $H_0^1(\Omega)$;
\item[(M5)] there exists a curve $\gamma_0(t)\in \Gamma$ passing through $u_0$ at $t=t_0$ and satisfying
    $$I_0(u_0)>I_0(\gamma_0(t))\;\hbox{for all\;}t\neq t_0.$$
\end{itemize}
Similar to \cite{JeongSeok.2013,CeramiZhongZou.2014}, we define a modified mountain pass energy level of $I_\lambda$
$$c_\lambda:=\min_{\gamma\in \Gamma_M}\max_{0\leq t\leq 1}I_\lambda (\gamma(t)),$$
where $$\Gamma_M=\Big\{\gamma\in \Gamma|\sup_{0\leq t\leq 1}\|\gamma(t)\|\leq M\Big\}\;\hbox{with}\;M:=2\max\Big\{\sup_{u\in \mathcal{S}}\|u\|, \sup_{t\in [0,1]}\|\gamma_0(t)\|\Big\}\;\hbox{fixed}.$$
Then it is easy to see that $\gamma_0\in \Gamma_M$ and thus
$c_0:=\min_{\gamma\in \Gamma_M}\max_{0\leq t\leq 1}I_0 (\gamma(t)).$
Similar to \cite{CeramiZhongZou.2014,JeongSeok.2013},   we  can easily prove the following results although the functional is different. Since the proofs are  analogous, we omit the details  and  refer the readers to \cite[Lemma 5.1-5.3, Proposition 5.1-5.3]{CeramiZhongZou.2014}.

\bl\lab{2014-4-9-l1} $c_\lambda\geq c_0\;\hbox{and}\;\lim_{\lambda\rightarrow 0}c_\lambda=c_0.$
\el

\bl\lab{2014-4-9-l2}
For $\forall\;d>0$, let $\{u_j\}\subset \mathcal{S}^d$, then up to a subsequence, $u_j\rightharpoonup u\in \mathcal{S}^{2d}$.
\el

\bl\lab{2014-4-9-l3}
Suppose that there exist sequences $\lambda_j<0, \lambda_j\rightarrow 0$ and $\{u_j\}\subset \mathcal{S}^d$ satisfying
$\underset{j\rightarrow \infty}{\lim}I_{\lambda_j}(u_j)\leq c_0\;\hbox{and}\;\underset{j\rightarrow \infty}{\lim}I'_{\lambda_j}(u_j)=0.$
Then there is $d_1>0$ such that for $0<d<d_1$, $\{u_j\}$ converges to some $u\in \mathcal{S}$ up to a subsequence.
\el

Next, we define
$m_\lambda:=\underset{0\leq t\leq 1}{\max}I_\lambda(\gamma_0(t)).$
Then by the definition of $c_\lambda$ we have that $c_\lambda\leq m_\lambda$. It is easy to see that $\underset{\lambda\rightarrow 0}{\lim}m_\lambda\leq c_0$. Combining with the conclusion of Lemma \ref{2014-4-9-l1}, we obtain that
\be\lab{2014-4-9-e1}
c_\lambda\leq m_\lambda\;\hbox{and}\;\underset{\lambda\rightarrow 0}{\lim}c_\lambda=\underset{\lambda\rightarrow 0}{\lim}m_\lambda=c_0.
\ee
We also define
$I_{\lambda}^{m_\lambda}:=\{u\in H_0^1(\Omega)|I_\lambda(u)\leq m_\lambda\}.$

\bl\lab{2014-4-9-l4}
For any $d_2,d_3>0$ satisfying $d_3<d_2<d_1$, there are constant $\alpha>0$ and $\lambda_0<0$ depending on $d_2,d_3$ such that for $\lambda\in (\lambda_0,0)$, we have
$$\|I'_\lambda(u)\|\geq \alpha\;\hbox{for all}\;u\in I_{\lambda}^{m_\lambda}\cap (\mathcal{S}^{d_2}\backslash \mathcal{S}^{d_3}).$$
\el

\bl\lab{2014-4-9-l5}
For $d>0$, there exists $\delta>0$ such that if $\lambda<0$ with $|\lambda|$ small enough,
$$t\in [0,1], I_\lambda(\gamma_0(t))\geq c_\lambda-\delta\;\hbox{implies}\;\gamma_0(t)\in \mathcal{S}^d.$$
\el

\bl\lab{2014-4-9-l6}
For any $d>0$ small and
$\lambda<0$ with $|\lambda|$ small enough depending on $d$, there exists a sequence $\{u_j\}\subset \mathcal{S}^d\cap I_{\lambda}^{m_\lambda}$ such that
$I'_\lambda(u_j)\rightarrow 0\;\hbox{as}\;j\rightarrow \infty.$
\el

\noindent{\bf Proof of Theorem\ref{main-th2}. } Taking $d>0$ small enough,
by Lemma \ref{2014-4-9-l6}, there exists some small $\lambda_0<0$ such that for   $\lambda_0<\lambda<0$, there exits a Palais-Smale sequence $\{u_j^\lambda\}\subset \mathcal{S}^{\frac{d}{2}}$. Since $\mathcal{S}$ is compact, it is easy to see $\{u_j^\lambda\}$ is bounded in $H_0^1(\Omega)$. Then by Lemma \ref{2014-4-9-l2}, there exists some $u^\lambda\in \mathcal{S}^{2\cdot \frac{d}{2}}=\mathcal{S}^{d}$ such that $u_j^\lambda\rightharpoonup u^\lambda$ up to a subsequence. Then we obtain $I'_\lambda(u^\lambda)=0$ and $u^\lambda\neq 0$. Hence $u^\lambda$ is a nontrivial critical point of $I_\lambda$. Testing by $u_-^\lambda$, we obtain that $\|u_-^\lambda\|^2=0$ which implies $u_-^\lambda\equiv 0$ and   $u^\lambda\geq 0$. Finally, by the maximum principle, we  know  that $u^\lambda>0$.
Hence,  $u^\lambda$ is a solution to (\ref{P1}).\hfill$\Box$

\end{CJK*}
 \end{document}